\documentclass{gtart_h}

\def\ifplaintex{\expandafter\ifx\csname documentclass\endcsname\relax}

\def\gtp{{\mathsurround=0pt\it $\cal G\mskip-2mu$eometry \&\ 
$\cal T\!\!$opology $\cal P\!$ublications}}  

\def\recd{{\small Received:\qua\receiveddate\ifx\reviseddate\relax
\else\qquad Revised:\qua\reviseddate\fi\par}} 

\def\lognumber#1{\def\thelognumber{#1}}
\def\volumenumber#1{\def\thevolumenumber{#1}}
\def\volumeyear#1{\def\thevolumeyear{#1}}
\def\papernumber#1{\def\thepapernumber{#1}}
\def\pagenumbers#1#2{\def\startpage{#1}\def\finishpage{#2}}
\def\published#1{\def\publishdate{#1}}

\def\received#1{\def\receiveddate{#1}}
\def\revised#1{\def\reviseddate{#1}}
\def\accepted#1{\def\accepteddate{#1}}

\long\def\asciiabstract#1{\long\def\theasciiabstract{#1}}

\let\\\par\let\thelognumber\relax\let\thevolumenumber\relax
\let\thepapernumber\relax\let\thevolumeyear\relax\let\startpage\relax
\let\finishpage\relax\let\publishdate\relax\let\receiveddate\relax
\let\reviseddate\relax\let\accepteddate\relax\let\theasciititle\relax
\let\theasciiauthors\relax
\let\theasciiabstract\relax

\let\theasciiemail\relax

\ifplaintex
\font\logobig=cmssbx10 scaled 3836
\font\logomed=cmssbx10 scaled 2557
\else
\font\logobig=cmssbx10 scaled 4200
\font\logomed=cmssbx10 scaled 2800
\fi

\long\def\makeagttitle{   
\count0=\startpage
\agt\hfill      
\hbox to 45truept{\vbox to 0pt{\vglue -13truept{\logomed A\kern -.37em{\logobig 
T}\kern -.38em G}\vss}\hss}
\break
{\small Volume \thevolumenumber\ (\thevolumeyear)
\startpage--\finishpage\nl
Published: \publishdate}

\vglue .25truein

{\parskip=0pt\leftskip 0pt plus
1fil\def\\{\par\smallskip}{\Large\bf\thetitle}\par\medskip} \vglue
0.05truein

{\parskip=0pt\leftskip 0pt plus 1fil\def\\{\par}{\sc\theauthors}
\par\medskip}%
 
\vglue 0.03truein 

{\small\leftskip 25truept\rightskip 25truept{\bf Abstract}\stdspace\theabstract

{\bf AMS Classification}\stdspace\theprimaryclass
\ifx\thesecondaryclass\relax\else; \thesecondaryclass\fi\par
{\bf Keywords}\stdspace \thekeywords\par}\vglue 7truept

}   

\ifplaintex
\font\phead=cmsl9 scaled 950
\font\pnum=cmbx10 scaled 913
\font\pfoot=cmsl9 scaled 950
\headline{\vbox to 0pt{\vskip -4.5mm\line{\small\phead\ifnum
\count0=\startpage ISSN 1472-2739 (on-line) 1472-2747 (printed)
\hfill {\pnum\folio}\else\ifodd\count0\def\\{ }%
\ifx\theshorttitle\relax\thetitle\else\theshorttitle\fi\hfill{\pnum\folio}
\else\def\\{ and }{\pnum\folio}\hfill\ifx\theshortauthors\relax\theauthors
\else\theshortauthors\fi\fi\fi}\vss}}
\footline{\vbox to 0pt{\vglue 0mm\line{\small\pfoot\ifnum\count0=\startpage
\copyright\ \gtp\hfill\else
\agt, Volume \thevolumenumber\ (\thevolumeyear)\hfill\fi}\vss}}
\else
\font=cmsl9 scaled 1050
\font\lnum=cmbx10 
\font=cmsl9 scaled 1050
\makeatletter
\def\@oddhead{{\small\lhead\ifnum\count0=\startpage ISSN 1472-2739 
(on-line) 1472-2747 (printed)\hfill {\lnum\number\count0}\else\ifodd\count0
\def\\{ }\ifx\theshorttitle\relax \thetitle \else\theshorttitle\fi\hfill
{\lnum\number\count0}\else\def\\{ and }{\lnum\number\count0}
\hfill\ifx\theshortauthors\relax 
\theauthors\else\theshortauthors\fi\fi\fi}}\def\@evenhead{\@oddhead}
\def\@oddfoot{\small\lfoot\ifnum\count0=\startpage\copyright\ \gtp\hfill\else
\agt, Volume \thevolumenumber\ (\thevolumeyear)\hfill\fi}
\def\@evenfoot{\@oddfoot}
\makeatother
\fi
\let\maketitlepage\makeagttitle
\let\makeshorttitle\maketitlepage
\let\maketitle\maketitlepage

\newwrite\gtoutfile
\long\gdef\makeheadfile{  
{\def\\{, }\def\s{ }
\immediate\openout\gtoutfile head.xxx
\immediate\write\gtoutfile{Proxy-for: \ifx\theasciiauthors\relax
\theauthors\else\theasciiauthors\fi\s<\ifx\theasciiemail\relax\theemail\else\theasciiemail\fi>}
\immediate\write\gtoutfile{\noexpand\\}
\immediate\write\gtoutfile{Authors: \ifx\theasciiauthors\relax
\theauthors\else\theasciiauthors\fi}
{\def\\{ }\immediate\write\gtoutfile{Title: \ifx\theasciititle\relax
\thetitle\else\theasciititle\fi}}
\immediate\write\gtoutfile{Subj-class: GT or SG, GR etc}
\immediate\write\gtoutfile{MSC-class: \theprimaryclass\ifx\thesecondaryclass\relax\else, \thesecondaryclass\fi}
\immediate\write\gtoutfile{Journal-ref: Algebr. Geom. Topol. \thevolumenumber\s
(\thevolumeyear) \startpage-\finishpage}
\immediate\write\gtoutfile{Comments: Published by Algebraic and
Geometric Topology at}
\immediate\write\gtoutfile{\s\s\s  http://www.maths.warwick.ac.uk/agt/AGTVol\thevolumenumber/agt-\thevolumenumber-\thepapernumber.abs.html}
\immediate\write\gtoutfile{\noexpand\\}
\immediate\write\gtoutfile{}
\ifx\theasciiabstract\relax
\immediate\write\gtoutfile{\theabstract}\else
\immediate\write\gtoutfile{\theasciiabstract}\fi
\immediate\write\gtoutfile{}
\immediate\write\gtoutfile{\noexpand\\}
\immediate\write\gtoutfile{}
\immediate\closeout\gtoutfile}}  

\def\maketitlepage{\makeagttitle\makeheadfile}
\let\makeshorttitle\maketitlepage
\let\maketitle\maketitlepage

\lognumber{32}
\volumenumber{4}
\volumeyear{2004}
\papernumber{32}
\pagenumbers{721}{755}
\received{2 April 2004} 
\revised{24 August 04}
\accepted{3 September 2004}
\published{11 September 2004}

\usepackage{amsmath, amssymb}
\usepackage[mathscr]{eucal}

\theoremstyle{plain}
\newtheorem{thm}{Theorem}[section]
\newtheorem{cor}[thm]{Corollary}
\newtheorem{prop}[thm]{Proposition}
\newtheorem{lemma}[thm]{Lemma}
\newtheorem{conj}[thm]{Conjecture}

\theoremstyle{definition}
\newtheorem*{rem}{Remark}

\newenvironment{smallgroup}
{\footnotesize}{\relax}

\DeclareMathOperator*{\Aut}{Aut} 
\DeclareMathOperator*{\Inn}{Inn}
\DeclareMathOperator{\Out}{Out}
\DeclareMathOperator*{\Stab}{Stab}

\DeclareMathOperator*{\Ima}{Im} 
\DeclareMathOperator*{\RE}{Re}
 
\DeclareMathOperator*{\tr}{Tr}

\DeclareMathOperator*{\GL}{GL}

\DeclareMathOperator*{\Ort}{O}

\DeclareMathOperator*{\SO}{SO}

\DeclareMathOperator*{\Uni}{U}
\DeclareMathOperator*{\PU}{PU}

\DeclareMathOperator*{\Sp}{Sp}

\DeclareMathOperator*{\Aff}{Aff}

\DeclareMathOperator*{\res}{Res}
\DeclareMathOperator*{\Isom}{Isom}

\newcommand{\eps}{\varepsilon}
\newcommand{\vp}{\varphi}

\newcommand{\al}{\alpha}

\newcommand{\be}{\beta}

\newcommand{\ga}{\gamma}
\newcommand{\Ga}{\Gamma}

\newcommand{\te}{\theta}

\newcommand{\si}{\sigma}

\newcommand{\om}{\omega}

\newcommand{\La}{\Lambda}

\newcommand{\ol}{\overline}
\newcommand{\wt}{\widetilde}
\newcommand{\wh}{\widehat}

\newcommand{\iny}{\infty}
\newcommand{\prt}{\partial}

\newcommand{\tri}{\ensuremath{\triangle}}

\newcommand{\innp}[1]{\left< #1 \right>}
\newcommand{\abs}[1]{\left\vert#1\right\vert}
\newcommand{\set}[1]{\left\{#1\right\}}
\newcommand{\brac}[1]{\left[#1\right]}
\newcommand{\pr}[1]{\left( #1 \right) }
\newcommand{\norm}[1]{\left\Vert#1\right\Vert}

\newcommand{\su}{\subset}

\newcommand{\op}{\oplus}

\newcommand{\smin}{\setminus}

\newcommand{\bdef}{\overset{\text{def}}{=}}

\newcommand{\lra}{\longrightarrow}
\newcommand{\lmto}{\longmapsto}

\newcommand{\B}[1]{\ensuremath{\mathbf{#1}}}
\newcommand{\BB}[1]{\ensuremath{\mathbb{#1}}}
\newcommand{\Cal}[1]{\ensuremath{\mathcal{#1}}}
\newcommand{\Fr}[1]{\ensuremath{\mathfrak{#1}}}

\newcommand{\Hy}{\ensuremath{\B{H}}}
\newcommand{\N}{\ensuremath{\BB{N}}}
\newcommand{\Q}{\ensuremath{\BB{Q}}}
\newcommand{\R}{\ensuremath{\B{R}}}
\newcommand{\Z}{\ensuremath{\BB{Z}}}
\newcommand{\C}{\ensuremath{\B{C}}}

\newcommand{\refP}[1]{Proposition~\ref{P:#1}}

\newcommand{\refT}[1]{Theorem~\ref{T:#1}}
\newcommand{\refL}[1]{Lemma~\ref{L:#1}}

\newcommand{\refC}[1]{Corollary~\ref{C:#1}}
\newcommand{\refE}[1]{Equation~(\ref{E:#1})}

\newcommand{\refCo}[1]{Conjecture~\ref{Co:#1}}

\begin{document}


\title{Peripheral separability and cusps of\\arithmetic 
hyperbolic orbifolds}
\author{D.B. McReynolds}
\address{University of Texas, Austin, TX 78712, USA}
\email{dmcreyn@math.utexas.edu}


\begin{abstract}
For $X = \R$, $\C$, or $\B{H}$, it is well known that cusp
cross-sections of finite volume $X$--hyperbolic $(n+1)$--orbifolds are
flat $n$--orbifolds or almost flat orbifolds modelled on the 
$(2n+1)$--dimensional Heisenberg group $\Fr{N}_{2n+1}$ or the 
$(4n+3)$--dimensional quaternionic Heisenberg group 
$\Fr{N}_{4n+3}(\B{H})$. We give a necessary and sufficient
condition for such manifolds to be diffeomorphic to a cusp
cross-section of an arithmetic $X$--hyperbolic $(n+1)$--orbifold.

A principal tool in the proof of this classification theorem is a subgroup
separability result which may be of independent interest.
\end{abstract}
\asciiabstract{%
For X = R, C, or H it is well known that cusp cross-sections of finite
volume X-hyperbolic (n+1)-orbifolds are flat n-orbifolds or almost
flat orbifolds modelled on the (2n+1)-dimensional Heisenberg group
N_{2n+1} or the (4n+3)-dimensional quaternionic Heisenberg group
N_{4n+3}(H). We give a necessary and sufficient condition for such
manifolds to be diffeomorphic to a cusp cross-section of an arithmetic
X-hyperbolic (n+1)-orbifold.  A principal tool in the proof of this
classification theorem is a subgroup separability result which may be
of independent interest.}




\maketitle


\section{Introduction}

\subsection{Main results}

A classical question in topology is whether a compact manifold
bounds. Hamrick and Royster \cite{HamrickRoyster82} showed
that every flat manifold bounds, and it was conjectured in
\cite{FarrellZdravkovska83} that any almost flat manifold bounds (for
some progress on this see see \cite{McReynolds04C} and \cite{Upadhyay01}). In
\cite{FarrellZdravkovska83}, Farrell and Zdravkovska made a stronger
geometric conjecture:  

\begin{conj}\label{Co:11}$\phantom{9}$
\begin{itemize}
\item[\rm(a)]
If $M^n$ is a flat Riemannian manifold, then $M^n =\prt W^{n+1}$ where
$W \smin \prt W$ supports a complete hyperbolic structure with finite
volume. 
\item[\rm(b)]
If $M^n$ supports an almost flat structure, then $M^n = \prt W^{n+1}$,
where $W \smin \prt W$ supports a complete Riemannian metric with
finite volume of whose sectional curvatures are negative. 
\end{itemize} 
\end{conj}

We say that a flat manifold $M^n$ \emph{geometrically bounds} if
(a) in \refCo{11} holds. Long and Reid \cite{LongReid00} showed that (a) is
false by proving that for a flat $(4n-1)$--manifold to geometrically bound the 
$\eta$--invariant is an integer. Furthermore, flat 3--manifolds with
nonintegral $\eta$--invariant are easily constructed using
\cite{Ouyang94}. Equivalently, this result of Long and Reid shows
that some flat manifolds cannot be diffeomorphic to a cusp
cross-section of a 1--cusped, finite volume real hyperbolic
4--manifold. On the other hand, Long and Reid showed \cite{LongReid02} that
every flat $n$--manifold is diffeomorphic to a cusp cross-section of
an arithmetic real hyperbolic $(n+1)$--orbifold.

For $X = \R$, $\C$, or $\B{H}$, cusp cross-sections of finite volume
$X$--hyperbolic $(n+1)$--orbifolds are flat $n$--manifolds or almost
flat orbifolds modelled on the $(2n+1)$--dimensional Heisenberg group
$\Fr{N}_{2n+1}$ or the $(4n+3)$--dimensional quaternionic Heisenberg
group $\Fr{N}_{4n+3}(\B{H})$. The first main result of this article 
shows that the result of Long and Reid in \cite{LongReid02} does
not generalize to the complex or quaternionic settings. Namely, (see
\S 2 for definitions):
 
\begin{thm}\label{T:Ob}$\phantom{9}$
\begin{itemize}
\item[\rm(a)]
For every $n \geq 2$, there exist infinite families of closed almost
flat $(2n+1)$--manifolds modelled on $\Fr{N}_{2n+1}$ which are not
diffeomorphic to a cusp cross-section of any arithmetic complex
hyperbolic $(n+1)$--orbifold. 
\item[\rm(b)]
For every $n \geq 1$, there exist infinite families of closed almost
flat $(4n+3)$--manifolds modelled on $\Fr{N}_{4n+3}(\B{H})$ which are
not diffeomorphic to a cusp cross-section of any finite volume
quaternionic hyperbolic $(n+1)$--orbifold. 
\end{itemize}
\end{thm}

Since all lattices in the isometry group of quaternionic hyperbolic
space are arithmetic (see \cite{Corlette92}), we drop the arithmeticity
assumption in (b). 

In order to give a complete classification of cusp cross-sections of
arithmetic hyperbolic lattices, we require certain subgroup
separability results. Recall that if $G$ is a group, $H<G$ and $g \in
G \smin H$, we say $H$ and $g$ are \emph{separated} if there exists a
subgroup $K$ of finite index in $G$ which contains $H$ but not $g$. We
say that $H$ is \emph{separable} in $G$ or $G$ is
\emph{$H$--separable}, if every $g \in G \smin H$ and $H$ can be
separated. We say that $G$ is \emph{LERF} (locally extendable
residually finite) if every finitely generated subgroup is separable.

We defer the statement of our second main result until \S 3 (see
\refT{Borel}) as it requires the language of algebraic
groups. Instead we state the result specialized to the rank--1
setting. For the statement, let $Y=\Hy_\R^n$, $\Hy_\C^n$,
$\Hy_\B{H}^n$ or $\Hy_\B{O}^2$.  

\begin{thm}[Stabilizer subgroup separability theorem]\label{T:Sep}
Let $\La$ be an arithmetic lattice in $\Isom(Y)$ and $v \in
\prt Y$. Then every subgroup of $\La \cap \Stab(v)$ is separable
in $\La$.    
\end{thm}

One well known application of subgroup separability is the lifting of
an immersion to an embedding in a finite cover (see \cite{Long87},
\cite{Hamilton01}, or \cite[p.\ 176]{MaclachlanReid03}). In the rank--1
setting, we have (see \refT{GImmersion} for a more general result):  

\begin{thm}\label{T:Immersion}
Let $\rho\co N \lra M$ be a $\pi_1$--injective immersion of an almost
flat manifold $N$ modelled in $\Fr{N}_{\ell n-1}(X)$ into an
arithmetic $X$--hyperbolic $m$--orbifold. Then there exists a finite
cover $\psi\co M^\prime \lra M$ such that $\rho$ lifts to an embedding.
\end{thm} 

A geometric corollary of \refT{Sep} of particular interest to us is: 

\begin{thm}\label{T:Cusp}$\phantom{9}$
\begin{itemize}
\item[\rm(a)]
A flat $n$--manifold is diffeomorphic to a cusp cross-section of an
arithmetic real hyperbolic $(n+1)$--orbifold if and only if
$\pi_1(M^n)$ injects into an arithmetic real hyperbolic
$(n+1)$--lattice.
\item[\rm(b)]
An almost flat $(2n+1)$--manifold  $M^{2n+1}$ modelled on
$\Fr{N}_{2n+1}$ is diffeomorphic to a cusp cross-section of an
arithmetic complex hyperbolic $(n+1)$--orbifold if and only if
$\pi_1(M^{2n+1})$ injects into an arithmetic complex hyperbolic
$(n+1)$--lattice.
\item[\rm(c)]
An almost flat $(4n+3)$--manifold $M^{4n+3}$ modelled on
$\Fr{N}_{4n+3}(\B{H})$ is diffeomorphic to a cusp cross-section of a
quaternionic hyperbolic $(n+1)$--orbifold if and only if
$\pi_1(M^{2n+1})$ injects into a quaternionic hyperbolic
$(n+1)$--lattice. 
\item[\rm(d)]
An almost flat $15$--manifold $M^{15}$ modelled on
$\Fr{N}_{15}(\B{O})$ is diffeomorphic to a cusp cross-section of an
octonionic hyperbolic $16$--orbifold if and only if
$\pi_1(M^{15})$ injects into an octonionic hyperbolic $16$--lattice. 
\end{itemize}
\end{thm} 

\refT{Cusp} reduces the classification of cusp-cross sections of
arithmetic $X$--hyperbolic $n$--orbifolds to the construction of faithful
representations of almost flat manifold groups into lattices. We
postpone stating the classification until \S 5 (see \refT{324}) as it
requires additional terminology. However, one interesting special case
which we state here is (see \S 7 for a proof): 

\begin{cor}\label{C:Nil3}
Every nil 3--manifold is diffeomorphic to a cusp cross-section of an
arithmetic complex hyperbolic 2--orbifold.
\end{cor}

The rest of the paper is organized as follows. We establish notation
and collect the preliminary material in \S 2 needed in the
sequel. Our main separability results are established in \S 3 and
\S 4 along with some algebraic and geometric corollaries. In \S 5 we
classify cusp cross-section of arithmetic hyperbolic orbifolds and
give the families of \refT{Ob} in \S 6. We conclude this article with
a detailed treatment of the nil 3--manifold case in \S 7.

\subsection{Acknowledgments}

I would like to thank my advisor Alan Reid for all his help. I am
indebted to Daniel Allcock for several helpful suggestions, most
notably the use of central products.  In addition, I would like to thank Karel
Dekimpe, Yoshinobu Kamishima, Richard Kent, and Richard Schwartz for conversations on this
work. Finally, I would like to thank the referee for several valuable
comments and for informing me of \refP{3}.

\section{Preliminary material}

In this section, we rapidly develop the material needed in the sequel.

\subsection{$X$--hyperbolic $n$--space}

For a general reference on this material, see \cite[II.10]{BridsonHaefliger99}. 
In all that follows, we let $X=\R,\C$, or $\B{H}$ and $\ell= \dim_\R X$.

Equip $X^{n+1}$ with a Hermitian form $H$ of signature $(n,1)$. We define
\emph{$X$--hyperbolic $n$--space} to be the (left) $X$--projectivization of the
$H$--negative vectors with the Bergmann metric associated to $H$. We
denote $X$--hyperbolic $n$--space together with this metric by
$\Hy_X^n$ and say that $\Hy_X^n$ is \emph{modelled on $H$} or call $H$
a \emph{model form}.

The \emph{boundary} of $\Hy_X^n$ in $PX^{n+1}$ is the
$X$--projectivization of the $H$--null vectors. We denote this set by
$\prt \Hy_X^n$, which is topologically just $S^{\ell n}$ (see 
\cite[p.\ 265]{BridsonHaefliger99}) and call the elements of the boundary
\emph{light-like vectors}. 

\subsection{The isometry group and lattices}

The isometry group of $\Hy_X^n$ is denoted by
$\Isom(\Hy_X^n)$. In each setting, $\Isom(\Hy_X^n)$ is locally
isomorphic to $\Uni(H)$. Specifically, 
\[ \Isom(\Hy_X^n) = \begin{cases} \innp{\PU(H)_0,\iota}, & X=\R,\C \\
\PU(H), & X=\B{H}, \end{cases} \]
where $\iota$ is an involution induced by inversion in the
real case and complex conjugation in the complex case. The usual
trichotomy for isometries holds in 
$\Isom(\Hy_X^n)$ (see \cite[p.\ 180--185]{Ratcliffe94}, \cite[p.\ 203]{Goldman99},
\cite{KimParker03}). Specifically, every (nontrivial) isometry is either
\emph{elliptic}, \emph{parabolic}, or \emph{loxodromic}.

We say that $\Ga < \Isom(\Hy_X^n)$ is a \emph{lattice} if $\Ga$ is a
discrete subgroup and $\Hy_X^n/\Ga$ has finite volume. In this case,
$M=\Hy_X^n/\Ga$ is called an \emph{$X$--hyperbolic
$n$--orbifold}. That finite volume manifolds exist, both compact and
noncompact, was established by Borel \cite{Borel63}. 

The spaces constructed in this way yield every locally
symmetric space of rank--1 except for those modelled on the
exceptional \emph{Cayley hyperbolic plane} $\Hy_\B{O}^2$. We shall
only make use of the fact that $\Isom(\Hy_\B{O}^2)$ has a faithful
linear representation and refer the reader to \cite{Allcock99} for more on
the Cayley hyperbolic plane. 

\subsection{The Heisenberg group and its quaternionic analog}

In the next few subsections, we introduce the $X$--Heisenberg group
and its automorphism group. See \cite{Goldman99} for the complex case
and \cite{KimParker03} for the quaternionic case. A thorough treatment
of this topic can be found in \cite{McReynolds04C}.

Let $\innp{\cdot,\cdot}$ denote the standard Hermitian product on $X^n$
and let $\om=\Ima \innp{\cdot,\cdot}$ be the associated hyper-symplectic 
form. The \emph{$X$--Heisenberg group 
$\Fr{N}_{\ell n-1}(X)$} is defined to be the topological space $X^{n-1} 
\times \Ima X$ together with the group structure
\[ (\xi_1,t_1)\cdot(\xi_2,t_2) \bdef
(\xi_1+\xi_2,t_1+t_2+2\om(\xi_1,\xi_2)). \]
The Lie group $\Fr{N}_{\ell n-1}(X)$ is simply connected and connected. 
Moreove, the group $\Fr{N}_{\ell n-1}(X)$ is nilpotent of step size two
in the case $X \ne \R$ and abelian in the case $X=\R$. The center (in
the nonabelian cases) is the commutator subgroup, which can be
identified with $\set{0} \times \Ima X$.

\subsection{Automorphisms of the $X$--Heisenberg group}

The automorphism group of the $X$--Heisenberg group $\Aut(\Fr{N}_{\ell
n-1}(X))$ splits as $\Inn(\Fr{N}_{\ell n-1}) \rtimes
\Out(\Fr{N}_{\ell n-1})$. The inner automorphism group can be identified 
with the vectors in $X^{n-1}$ which are $\om$--nondegenerate together
with the zero vector. In the real case, this set is just $\set{0}$,
while in the other two cases, this is the whole of $X^{n-1}$. 

The outer automorphism group is comprised of three types of
automorphisms. The first type of automorphism is a \emph{symplectic
rotation} given by $S(\xi,t) = (S\xi,t)$ for $S \in \Sp(\om)$. 
The second type of automorphism is a \emph{Heisenberg dilation} given by
$d(\xi,t) = (d\xi,d^2t)$ for $d \in \R^\times$. Finally, we have 
\emph{$X$--scalar conjugation} given by $\zeta(\xi,t) = (\zeta^{-1}\xi\zeta,
\zeta^{-1}t\zeta)$ for $\zeta \in X^\times$. The outer automorphism 
group is generated by these three automorphisms. In summary, we have
\[ \Out(\Fr{N}_{\ell n-1}(X)) = \begin{cases} \GL(n-1;\R), & X=\R \\
  \Sp(2n-2) \times \R^\times, & X=\C \\ \Sp(\om) \times \R^\times
  \times \B{H}^\times, & X=\B{H}. \end{cases} \]

\subsection{Maximal compact subgroups}

Our primary concern is with maximal compact subgroups
$\Aut(\Fr{N}_{\ell n-1})$. The maximal compact subgroups are of the form
\[ M(X) = \begin{cases} \Ort(B_{M(X)}), & X=\R \\ \innp{\Uni(H_{M(X)}),\iota}, &
  X=\C \\ \Uni(H_{M(X)}) \times S, & X=\B{H}, \end{cases} \]
where $B_{M(X)}$ is a symmetric, positive definite bilinear form,
$H_{M(X)}$ is a signature $(n-1,0)$ Hermitian form with $\Ima H_{M(X)}
= \om$, and $S$ is the unit sphere in $\B{H}$ (equipped possibly with
a nonstandard quaternionic structure). Since the maximal
compact subgroups are conjugate, each $M(X)$ is conjugate to
\[ M_s(X) = \begin{cases} \Ort(n-1), & X=\R \\ \innp{\Uni(n-1),\iota},
  & X=\C \\ \Sp(n-1) \times \SO(3), & X=\B{H}. \end{cases} \]

For a given maximal compact subgroup $M$, we call the group
$\Fr{N}_{\ell n-1}(X) \rtimes M$ a \emph{unitary affine group} and
denote this group by $U_M(n-1;X)$. We call the group $\Fr{N}_{\ell
n-1}(X) \rtimes (M(X) \times \R^+)$ an \emph{$X$--Heisenberg
similarity group} and denote this group by $S_M(n-1;X)$. Finally, we
call the group $\Fr{N}_{\ell n-1}(X) \rtimes \Aut(\Fr{N}_{\ell
  n-1}(X))$ the \emph{$X$--Heisenberg affine group} and denote it by
$\Aff(\Fr{N}_{\ell n-1}(X))$.  


\subsection{Almost crystallographic groups modelled on the
\newline $X$--Heisenberg group} 

In this subsection, we introduce almost crystallographic groups
modelled on the $X$--Heisenberg group. We refer the reader to
\cite{Dekimpe96} or \cite[Chapter II and VIII]{Raghunathan72} for a
general treatment on discrete subgroups in nilpotent Lie groups.

By an \emph{almost crystallographic group} or \emph{AC-group} modelled
on $\Fr{N}_{\ell n-1}(X)$, we mean a discrete subgroup $\Ga <
\Aff(\Fr{N}_{\ell n-1})(X)$ such that $\Fr{N}_{\ell
n-1}(X)/\Ga$ is compact and $\Ga \cap \Fr{N}_{\ell n-1}(X)$ is a
finite index subgroup of $\Ga$. When $\Ga$ is torsion free, we
say that $\Ga$ is an \emph{almost Bieberbach group} or \emph{AB-group}
modelled on $\Fr{N}_{\ell n-1}(X)$. Every AC-group modelled on
$\Fr{N}_{\ell n-1}(X)$ is determined by the short exact sequence 
\[ 1 \lra L \lra \Ga \lra \te \lra 1, \]
where $L = \Ga \cap \Fr{N}_{\ell n-1}(X)$ and $\abs{\te}<\iny$.
We call $L$ the \emph{Fitting subgroup of $\Ga$} and $\te$ the
\emph{holonomy group of $\Ga$}.

It is well known (see \cite[Ch. 3]{Dekimpe96}) that the above exact sequence
induces an injective homomorphism $\vp\co \te \lra \Out(\Fr{N}_{\ell
n-1}(X)) <  \Aut(\Fr{N}_{\ell n-1}(X))$ which we call the
\emph{holonomy representation of $\te$}. Since $\te$ is finite, this
is conjugate into a representation $\vp\co \te \lra M(X)$ for any
$M(X)$. This yields a faithful representation $\rho\co \Ga \lra
U_M(n-1;X)$ for any $M(X)$.

\subsection{Almost flat manifolds}

Let $(M^n,g)$ be a complete Riemannian manifold. We let $d=d(g)$,
$c^-(g)$ and $c^+(g)$ denote the diameter of $M$ and the lower and
upper bounds of the sectional curvature of $M$, respectively, and set 
$c(g)$ to be the maximum of $\abs{c^+}$ and $\abs{c^-}$. We say that
$M$ is \emph{almost flat} if there exists a family of complete
Riemannian metrics $g_j$ on $M$ such that 
\[ \lim_{j \lra \iny} d(g_j)^2c(g_j) = 0. \]
Gromov \cite{Gromov78} proved that every compact almost flat manifold
is of the form $N/\Ga$, where $N$ is a connected, simply connected
nilpotent Lie group and $\Ga$ is an AB-group modelled on $N$. 

Of importance to us is some of the generalized Bieberbach theorem (see
\cite{Dekimpe96}).

\eject
\begin{thm}[Generalized Bieberbach theorem]\label{T:GenBieberbach}$\phantom{9}$
\begin{itemize}
\item[\rm(a)]
Let $M$ be an almost flat manifold with universal cover $\Fr{N}_{\ell
  n-1}$. Then there exists a faithful representation
$\vp\co \pi_1(M) \lra \Aff(\Fr{N}_{\ell n-1}(X))$ such that
$\vp(\pi_1(M))$ is an AB-group. 
\item[\rm(b)]
$M=\Fr{N}_{\ell n-1}/\Ga$ and $M^\prime=\Fr{N}_{\ell n-1}/\Ga^\prime$
are diffeomorphic if and only if there exists $\al \in
\Aff(\Fr{N}_{\ell n-1})$ such that
\[ \Ga^\prime = \al^{-1}\Ga \al. \]
\end{itemize}
\end{thm}

In the remainder of this article we refer to compact almost flat
manifolds as \emph{infranil manifolds} modelled on $N$, where $N$ is
the connected, simply connected nilpotent cover. In the event the 
fundamental group is a lattice in $N$, we call such manifolds
\emph{nil manifolds} modelled on $N$.

\subsection{Maximal peripheral subgroups, stabilizer groups, and cusps}

For a lattice $\La < \Isom(\Hy_X^n)$ with cusp at $v$, we
define the \emph{maximal peripheral subgroup of $\La$ at $v$} to
be the subgroup $\tri_v(\La) = \Stab(v) \cap \La$.
This is the subgroup generated by the parabolic and elliptic
isometries of $\La$ fixing $v$. By the Kazhdan-Margulis theorem (this
is sometimes called Margulis' lemma; see \cite[Chapter
XI]{Raghunathan72}), $\tri_v(\La)$ is virtually nilpotent. Specifically,
the maximal nilpotent subgroup of $\tri_v(\La)$ is given by $L =
\tri_v(\La) \cap N$, where $N$ is isomorphic to $\Fr{N}_{\ell n-1}(X)$.
Moreover, the Kazhdan-Margulis theorem allows us to select a
horosphere $\Cal{H}$ such that $\Cal{H}/\tri_v(\La)$ is embedded in
$\Hy_X^n/\La$. In this case, we call $\Cal{H}/\tri_v(\La)$ a
\emph{cusp cross-section} of the cusp at $v$. Often when $v$ is
unimportant, we simply write $\tri(\La)$.

More generally, for any $v \in \prt \Hy_X^n$, we define $\tri_v(\La) =
\La \cap \Stab(v)$ and call this subgroup the \emph{stabilizer group
of $\La$ at $v$}. There are three possibilities: 
\begin{itemize}
\item[(1)]
$\tri_v(\La)$ is finite.
\item[(2)]
$\tri_v(\La)$ is virtually cyclic with cyclic subgroup generated by a
loxodromic isometry. 
\item[(3)]
$\tri_v(\La)$ is an AC-group modelled on $\Fr{N}_{\ell n-1}(X)$.
\end{itemize}

 
\subsection{Iwasawa decompositions of the isometry group}  

For the isometry group of $X$--hyperbolic $n$--space, we
can decompose $\Isom(\Hy_X^n)$ as $KAN$ via the \emph{Iwasawa
decomposition} (see \cite[p.\ 311--313]{BridsonHaefliger99}). The factor $N$
is isomorphic to the $X$--Heisenberg group $\Fr{N}_{\ell n-1}(X)$ and
all isomorphisms arise in the following fashion. Let $H$ be a model
Hermitian form for $X$--hyperbolic $n$--space and $V_\iny$ be the $H$--orthogonal 
complement of $v_0$ and $v_\iny$, a pair of $X$--linearly independent
$H$--null vectors in $X^{n+1}$. For a maximal compact group $M(X)$ with
associated Hermitian form $H_{M(X)}$, let
\[ \psi\co (X^{n-1},H_{M(X)}) \lra (V_\iny,H_{|V_\iny}) \]
be any isometric $X$--isomorphism. This induces a map
$\eta\colon X^{n-1} \lra N$ defined by $\eta(\xi) = \exp(\psi(\xi)
v_\iny^* - v_\iny\psi(\xi)^*)$, 
where $xy^*(\cdot) = H(\cdot,y)x$ is the Hermitian outer pairing of $x$ and $y$ 
with respect to the Hermitian form $H$. This extends to all of
$\Fr{N}_{\ell n-1}(X)$ as these elements generate $\Fr{N}_{\ell
  n-1}(X)$. In fact, this extends to $\eta\co S_M(n-1;X) \lra \Isom(\Hy_X^n)$.
Since these isometries preserve $v_\iny$, this yields
$\eta(S_M(n-1;X)) = \Stab(v_\iny)$. 

\subsection{Algebraic groups}

As we use the language of algebraic groups throughout this paper, in
this subsection we review some of the basic material. See
\cite{Borel91} or \cite[Ch. 2]{PlatonovRapinchuk94}.

In the remainder of this article, all fields are assumed to be
algebraic number fields unless stated otherwise.

By a \emph{linear algebraic group} we mean a subgroup of $\GL(n;\C)$
which is closed in the Zariski topology. We say that $G$ is
\emph{$k$--algebraic} when there is a generating set of
$k$--polynomials for $\Fr{a}_G$, the ideal vanishing on $G$. 
For any subring $R\su \C$, we define the \emph{$R$--points} of $G$ to be
the subgroup $G \cap \GL(n;R)$. We denote the $R$--points
of $G$ by $G_R$. 

A \emph{Borel subgroup} of $G$ is a maximal, connected solvable
subgroup of $G$. Borel subgroups of $G$ are conjugate in $G$ and
conjugate into the subgroup of upper triangular matrices. If $G$ is  
$k$--algebraic, then $B$ will be $k^\prime$--algebraic for some finite
extension $k^\prime$ of $k$.  

A \emph{maximal algebraic torus} $T$ of $G$ is a maximal
diagonalizable algebraic subgroup. If $k$ is the 
field of definition for $G$, then the \emph{splitting field}
$k^\prime$ for $T$ is a finite extension of $k$. This is the smallest
field for which $T$ can be diagonalized. In particular, $T$ will be a
$k^\prime$--algebraic group. We say that $U<G$ is \emph{unipotent} if
$U$ is conjugate to a subgroup of the upper triangular matrices with
ones along the diagonal. Maximal unipotent subgroups $U$ are
connected, nilpotent, algebraic subgroups and if $G$ is
$k$--algebraic, $U$ is $k^\prime$--algebraic for some finite extension
$k^\prime$ of $k$. We note that every maximal torus $T$ or maximal
unipotent subgroup $U$ in $G$ is contained in a Borel subgroup (see
\cite[Cor.\ 11.3]{Borel91}). 

Finally, we require the following lemma in the sequel and refer the
reader to \cite[Cor.\ 10.14]{Raghunathan72} for a proof. In the
statement, $\Cal{O}_k$ denotes the ring of algebraic integers in the
number field $k$ (see \cite{Weil95}).

\begin{lemma}\label{L:323}
Let $f\co G \lra G^\prime$ be a $k$--homomorphism of $k$--algebraic
groups. If $\Ga< G_k$ is commensurable with $G_{\Cal{O}_k}$, then
there exists $\Ga^\prime<G_k^\prime$, commensurable with
$G_{\Cal{O}_k}^\prime$ such that $f(\Ga) < \Ga^\prime$.
\end{lemma}

\section{Borel subgroup separability theorem}

This section is devoted to proving the following result. 

\begin{thm}[Borel subgroup separability theorem]\label{T:Borel}
Let $G$ be a connected $k$--algebraic group and $B$
a Borel subgroup of $G$. Then any subgroup of $B_{\Cal{O}_k}$ is
separable in $G_{\Cal{O}_k}$. 
\end{thm}

Before embarking upon the proof, we record some facts that will be
needed. We begin with the following lemma (see \cite{LongReid02}). 

\begin{lemma}\label{L:293}
Let $G$ be a group and $H<K<G$. If $H$ is separable in $G$ and
$[K:H]<\iny$, then $K$ is separable in $G$. 
\end{lemma}

\begin{lemma}\label{L:295}
Let $G$ be a group and assume that $H,L<G$ are separable in $G$. Then
$H \cap L$ is separable in $G$.
\end{lemma}

\begin{proof}
Let $\ga \in G \smin (H \cap L)$ and assume that $\ga \notin H$. Since
$H$ is separable in $G$, there exists a finite index subgroup $K<G$
with $H<K$ and $\ga \notin K$. As $H \cap L < H$, $K$ separates $\ga$
and $H \cap L$, as needed. For the alternative, an identical argument
is made. 
\end{proof}

\begin{lemma}\label{L:294}
Let $G$ be a group, $G_0$ a subgroup of finite index, and $H$ a
subgroup. $H$ is separable in $G$ if and only if $(G_0 \cap H)$ is
separable in $G_0$. 
\end{lemma}

\begin{proof}
The direct implication follows immediately from \refL{295},
since\break $H \cap G_0$ is separable in the larger group $G$. For the
reverse implication, to show that $H$ is separable in $G$, by
\refL{293} it suffices to show that $G_0 \cap H$ is separable in
$G$. For $g \in G \smin (G_0 \cap H)$, there are two cases to
consider. If $g \notin G_0$, then $G_0$ separates $G_0 \cap H$ and
$g$. Otherwise, if $g \in G_0$, since $G_0 \cap H$ is separable in
$G_0$, there exists a finite index subgroup $K<G_0$ such that $G_0
\cap H<K$ and $g \notin K$. Since $[G:G_0]<\iny$, $K$ is the desired
finite index subgroup of $G$ separating $G_0 \cap H$ and $g$.
\end{proof}

The separability of Borel subgroups relies on the following result of
Chahal \cite{Chahal80} which establishes the congruence subgroup property
for solvable algebraic groups defined over number fields. Before we state
the result, we recall the definition of congruence kernels and
reduction homomorphisms. 

For each ideal $\Fr{p}<\Cal{O}_k$, we can define the 
homomorphism (\emph{reduction modulo $\Fr{p}$}) 
$r_\Fr{p}\co H_{\Cal{O}_k} \lra \GL(m;\Cal{O}_k/\Fr{p})$
by $r_\Fr{p}(\ga) = (\ga_{ij} \mod \Fr{p})_{ij}$.
By a \emph{congruence kernel} we mean a subgroup $\ker r_\Fr{p}$, for
some (nontrivial) ideal $\Fr{p}<\Cal{O}_k$ and denote this subgroup by
$K_{H,\Fr{p}}$. 

\begin{thm}\label{T:Ch}
Let $H$ be a solvable $k$--algebraic group. Then every finite index
subgroup of $H_{\Cal{O}_k}$ contains a congruence kernel.
\end{thm}

\subsection{The proof of \refT{Borel}}

For the proof of \refT{Borel}, recall that $G$ is a connected
$k$--algebraic group with a Borel subgroup $B$ defined over
$k^\prime$. The strategy for the proof is as follows. If $B$ is
defined over $k$ ($G$ is $k$--split), the proof reduces to proving
that $B_{\Cal{O}_k}$ is separable in $G_{\Cal{O}_k}$. For once this
has been established, to separate a subgroup of $B_{\Cal{O}_k}$ in
$G_{\Cal{O}_k}$, it suffices to separate the subgroup in
$B_{\Cal{O}_k}$. The latter is achieved by appealing to a theorem of
Mal'cev. In the non-split case when $k^\prime$ is not contained in
$k$, we enlarge our field to the composite field of $k$ and $k^\prime$
and appeal to the split case. In the remainder of this subsection, we
give the details.

The following two lemmas comprise the key steps in the proof of
\refT{Borel}. 

\begin{lemma}\label{L:1}
Let $G$ be a connected $k$--algebraic group and $B$ a $k$--defined Borel
subgroup of $G$. If $B_{\Cal{O}_k}$ is separable in
$G_{\Cal{O}_k}$, then every subgroup of $B_{\Cal{O}_k}$ is separable
in $G_{\Cal{O}_k}$. 
\end{lemma}

\begin{lemma}\label{L:2}
Let $G$ be a connected $k$--algebraic group and $B$ a $k$--defined
Borel subgroup of $G$. Then $B_{\Cal{O}_k}$ is separable in
$G_{\Cal{O}_k}$.
\end{lemma} 

Assuming these lemmas, we prove \refT{Borel}.    

\begin{proof}[Proof of \refT{Borel}]
The proof breaks into two cases, depending on\break whether or not
$k^\prime\su k$.  

\medskip
\textbf{Case 1}\qua $k^\prime\su k$

Since $k^\prime \su k$, $B$ is a $k$--defined Borel subgroup
of $G$. Therefore by \refL{2}, $B_{\Cal{O}_k}$ is separable in
$G_{\Cal{O}_k}$. Thus by \refL{1}, every subgroup of $B_{\Cal{O}_k}$
is separable in $G_{\Cal{O}_k}$, as desired.

\medskip
\textbf{Case 2}\qua $k^\prime$ is not contained in $k$

In this case, let $\wh{k}$ denote the composite of $k$ and
$k^\prime$. Then $G$ is a $\wh{k}$--algebraic group and $B$ is a
$\wh{k}$--defined Borel subgroup. Therefore by \refL{2},
$B_{\Cal{O}_{\wh{k}}}$ is separable in $G_{\Cal{O}_{\wh{k}}}$. Thus by
\refL{1}, every subgroup of $B_{\Cal{O}_{\wh{k}}}$ is separable in
$G_{\Cal{O}_{\wh{k}}}$. Since $k \su \wh{k}$, $B_{\Cal{O}_k} \su
B_{\Cal{O}_{\wh{k}}}$ and so every subgroup of $B_{\Cal{O}_k}$ is
separable in $G_{\Cal{O}_{\wh{k}}}$. Thus every subgroup of
$B_{\Cal{O}_k}$ is separable in the smaller group $G_{\Cal{O}_k}$.
\end{proof} 
   
We are now left with the task of verifying \refL{1} and \refL{2}. 

\begin{proof}[Proof of \refL{1}]
Let $S<B_{\Cal{O}_k}$ be a subgroup. For $\ga \in G_{\Cal{O}_k} \smin
S$, there are two cases to consider. First, if $\ga \notin
B_{\Cal{O}_k}$, then by the separability of $B_{\Cal{O}_k}$ we can
find a finite index subgroup $K<G_{\Cal{O}_k}$ such that $S<B_{\Cal{O}_k}<K$
and $\ga \notin K$. If $\ga \in B_{\Cal{O}_k}$ we argue as
follows. Since $B_{\Cal{O}_k}$ is polycyclic (\cite[p.\ 53]{Raghunathan72} or
\cite[p.\ 196]{Wehrfritz73}) it is LERF by \cite{Malcev58}. Therefore there exists a
finite index subgroup $K_B<B_{\Cal{O}_k}$ such that $S<K_B$ and $\ga 
\notin K_B$. By \refT{Ch}, $B_{\Cal{O}_k}$ has the
congruence subgroup property. Thus there exists a congruence kernel
$K_{B,\Fr{p}}$ of $B_{\Cal{O}_k}$ with $K_{B,\Fr{p}}<K_B$. As
$K_{B,\Fr{p}}$ is the intersection of $B_{\Cal{O}_k}$ 
with the congruence kernel $K_{G,\Fr{p}}$ of $G_{\Cal{O}_k}$, by
\refL{295}, $K_{B,\Fr{p}}$ is separable in $G_{\Cal{O}_k}$. By
\refL{293}, $K_B$ is separable in $G_{\Cal{O}_k}$, since
$[K_B:K_{B,\Fr{p}}]<\iny$. Consequently, we can find a finite index
subgroup $K<G_{\Cal{O}_k}$ such that $S<K_B<K$ and $\ga \notin
K$. Therefore $S$ and $\ga$ are separated in $G_{\Cal{O}_k}$.
\end{proof}

The proof of \refL{2} follows from a more general result established
in \cite{Bergeron00} (see also \cite{MargulisSoifer81}):

\begin{prop}\label{P:3}
Let $H$ be an algebraic group in a linear algebraic group $G$ and
$\Ga$ a finitely generated subgroup of $G$. Then $H \cap \Ga$ is
separable in $\Ga$.
\end{prop}

\refL{2} follows from \refP{3} by setting $H=B$ and
$\Ga=G_{\Cal{O}_k}$.

The proof of \refT{Borel} works in greater generality. Specifically, 

\begin{cor}\label{C:Gen}
Let $G$ be a connected $k$--algebraic group and $N$ a $k^\prime$--algebraic
subgroup with $k \su k^\prime$. If $N_{\Cal{O}_{k^\prime}}$
has the congruence subgroup property, then a finitely
generated subgroup $L$ of $N_{\Cal{O}_k}$ is separable in
$G_{\Cal{O}_k}$ if and only if $L$ is separable in $N_{\Cal{O}_k}$.  
\end{cor}

\subsection{Corollaries to \refT{Borel}}

In this subsection, we state a few corollaries to \refT{Borel} pertaining
to general algebraic groups.

Our first corollary shows that the conclusions of \refT{Borel} hold for
any subgroup of $G$ commensurable with $G_{\Cal{O}_k}$. We call such
subgroups \emph{$k$--arithmetic subgroups}.

\begin{cor}\label{C:CBorel}
Let $G$ be a connected $k$--algebraic group, $\La$ a
$k$--arithmetic subgroup in $G$, and $B$ a Borel subgroup of $G$.
Then every subgroup of $\La \cap B$ is separable in $\La$.
\end{cor}

\begin{proof}
For a subgroup $S<B\cap \La$, by \refL{294}, it
suffices to separate $S \cap G_{\Cal{O}_k}$ in $G_{\Cal{O}_k} \cap
\La$. Since $S \cap G_{\Cal{O}_k}$ is a subgroup of $B_{\Cal{O}_k}$,
by \refT{Borel}, $S \cap G_{\Cal{O}_k}$ is separable in
$G_{\Cal{O}_k}$. Thus, $S \cap G_{\Cal{O}_k}$ is separable in
$G_{\Cal{O}_k} \cap \La$.
\end{proof}  

As a result of \refC{CBorel}, every corollary and theorem stated below
implies the same result for any $k$--arithmetic subgroup in
$G$. Consequently, we only state the results for group of
$k$--integral points. The connected assumption is unnecessary since
every $k$--algebraic group has finitely many connected components (see
\cite[p.\ 51]{PlatonovRapinchuk94}).   

One corollary to \refT{Borel} is:

\begin{cor}\label{C:3Cor}
Let $G$ be a connected $k$--algebraic group.
\begin{itemize}
\item[\rm(a)]
If $U<G$ is a maximal unipotent subgroup, then every subgroup of
$U_{\Cal{O}_k}$ is separable in $G_{\Cal{O}_k}$. 
\item[\rm(b)]
If $T<G$ is a maximal torus, then every subgroup of $T_{\Cal{O}_k}$ is 
separable in $G_{\Cal{O}_k}$.
\item[\rm(c)]
If $S<G_{\Cal{O}_k}$ is a solvable subgroup, then $S$ is separable in
$G_{\Cal{O}_k}$. 
\end{itemize}
\end{cor}

\begin{proof}
(a) and (b) follow immediately from \refT{Borel} since $U$ and $T$ are
contained in a Borel subgroup. For (c), since every solvable subgroup is virtually
contained in a Borel subgroup (see \cite[p.\ 137]{Borel91}), by \refL{293},
it suffices to separate $S \cap B$ in $G_{\Cal{O}_k}$. The latter is
done using \refT{Borel}.  
\end{proof}

For a $k$--algebraic group $G$, by an \emph{arithmetic $G$--orbifold},
we mean a topological manifold of the form $G/\La$, where $\La$ is an
arithmetic lattice in $G$.

\begin{thm}\label{T:GImmersion}
Let $\rho\co N \lra M$ be a $\pi_1$--injective immersion of an
infrasolv manifold $N$ into an arithmetic $G$--orbifold $M$. Then
there exists a finite cover $\psi\co M^\prime \lra M$ such that $\rho$
lifts to an embedding. 
\end{thm}

\begin{proof}
The map $\rho$ induces a homomorphism $\rho_*\co \pi_1(N) \lra
\pi_1(M)$. Since $N$ is an infrasolv manifold, $\rho_*(\pi_1(N))$ is a
solvable subgroup of $\pi_1(M)$. Since $\pi_1(M) = \La$, for some
arithmetic lattice in $G$, by \refC{3Cor}, $\rho_*(\pi_1(M))$ is
separable in $\pi_1(M)$. It now follows by a standard argument (see
\cite{Long87}) that $\rho$ can be promoted to an embedding in some
finite covering of $M$. 
\end{proof}


\section{The stabilizer subgroup separability theorem} 

In this section we prove \refT{Sep} and corollaries specific
to lattices in the isometry group of hyperbolic space. 

\subsection{Stabilizer subgroup separability}

As mentioned in \S 2.7, there is a simple trichotomy for the stabilizer
groups of light-like vectors for $X$--hyperbolic lattices. For a
lattice $\La < \Isom(\Hy_X^n)$ and $v \in \prt \Hy_X^n$, exactly one of 
the following holds:
\begin{itemize}
\item[(1)]
$\tri_v(\La)$ is finite.
\item[(2)]
$\tri_v(\La)$ is virtually cyclic with maximal cyclic subgroup
generated by a loxodromic isometry.
\item[(3)]
$\tri_v(\La)$ is an AC-group modelled on the $X$--Heisenberg group
$\Fr{N}_{\ell n -1}(X)$. 
\end{itemize}

\begin{proof}[Proof of \refT{Sep}]
To prove \refT{Sep}, we split our consideration naturally into three 
cases depending on the above trichotomy. 

Since $X$--hyperbolic lattices are residually finite it follows easily from
\refL{293} that subgroups in case (1) are separable. For $X=\R$ or
$\C$, case (2) follows exactly the proof in \cite{Hamilton01} on noting $\GL(n;\C) \lra
\GL(2n;\R)$. For $X=\B{H}$ or $\B{O}$, since every lattice in
$\Isom(\Hy_\B{H}^n)$ and $\Isom(\Hy_\B{O}^2)$ is arithmetic, we can
apply \refC{3Cor} (c) to separate. 

For (3), as peripheral subgroups are virtually unipotent, \refC{3Cor} 
handles this case. To be complete, we first realize the arithmetic
lattice $\La$ as a subgroup of $\GL(m;\Q)$ with a finite index
subgroup in $\GL(m;\Z)$ and finish by applying \refC{CBorel} with
\refC{3Cor}. 
\end{proof}

\begin{rem}
In \cite{Hamilton01}, Hamilton proved that in a cocompact lattice $\La <
\Isom(\Hy_\R^n)$, every virtually abelian subgroup is
separable. As her proof does not require arithmeticity, our
proof of \refT{Sep} uses arithmeticity only in (3). 
\end{rem}

\begin{cor}\label{C:RealAbel} 
Let $\La$ be an arithmetic real hyperbolic lattice and $A$ an abelian
subgroup. Then $A$ is separable in $\La$.
\end{cor} 

The analog of abelian subgroups in the complex, quaternionic,
octonionic settings are nilpotent subgroups. In the complex setting,
we have: 

\begin{cor}\label{C:ComplexNil}
Let $\La$ be an arithmetic complex hyperbolic lattice and $N$ a
nilpotent subgroup. Then $N$ is separable in $\La$.
\end{cor}

Since all lattices in $\Isom(\Hy_\B{H}^n)$ and $\Isom(\Hy_\B{O}^2)$ are 
arithmetic, we may drop the arithmeticity condition to obtain: 

\begin{cor}\label{C:QuatNil}
Let $\La$ be a lattice in $\Isom(\Hy_\B{H}^n)$ or
$\Isom(\Hy_\B{O}^2)$ and $N$ a nilpotent subgroup. Then $N$ is
separable in $\La$. 
\end{cor}  

\section{A necessary and sufficient condition for arithmetic
admissibility}

The goal of this section is to give a classification of cusp
cross-sections of arithmetic $X$--hyperbolic $n$--orbifolds. By
\refT{Cusp}, we are reduced to classifying AB-groups which admit injections into
arithmetic $X$--hyperbolic lattices. The main point of this section is
to prove that this is equivalent to constructing injections into
arithmetically defined subgroups of unitary affine groups. The latter
groups are easier to work with in regard to this problem, as the
generalized Bieberbach theorems ensure the existence of
injections. The proof of this reduction relies on being able to
realize unitary affine groups as algebraic subgroups in the isometry
group of $X$--hyperbolic space.  In total,
this section is straightforward with the bulk of the material consisting
of terminology, notation, and formal manipulation. We hope the main
point of this section is not lost in this.   
 
\subsection{Characterization of noncocompact arithmetic lattices}

In this subsection, we give the classification of noncocompact arithmetic
$X$--hyperbolic $n$--lattices. This is originally due to Weil \cite{Weil60}.
We refer the reader to \cite{McReynolds04C} for a proof.

\begin{thm}\label{T:Arithmetic}
Let $\La$ be a noncocompact arithmetic lattice in $\Isom(\Hy_X^n)$.
\begin{itemize}
\item[\rm(a)]
If $X=\R$, then $\La$ is conjugate to an arithmetic lattice in $\Ort(B)$, 
where $B$ is a signature $(n,1)$ bilinear form defined over $\Q$.
\item[\rm(b)]
If $X=\C$, then $\La$ is conjugate to an arithmetic lattice in 
$\Uni(H)$, where $H$ is a Hermitian form of signature $(n,1)$ defined 
over an imaginary quadratic number field.
\item[\rm(c)]
If $X=\B{H}$, then $\La$ is conjugate to an arithmetic lattice in 
$\Uni(H)$, where $H$ is a Hermitian form of signature $(n,1)$ defined
over a definite quaternion algebra with Hilbert symbol
$\pr{\frac{-a,-b}{\Q}}$ for $a,b\in\N$. 
\end{itemize}
\end{thm}

\subsection{Algebraic structure of unitary affine groups}

Recall for each maximal compact subgroup $M(X)$ of $\Aut(\Fr{N}_{\ell n-1})$,
we defined the unitary affine group $U_M(n-1;X)$ to be $\Fr{N}_{\ell
n-1}(X) \rtimes M(X)$. The algebraic structure of these groups is
completely determined by the algebraic structure of the maximal
compact subgroup. Specifically, $U_M(n-1;X)$ is $k$--algebraic if and
only if $M$ is $k$--algebraic. In turn, the algebraic structure of $M$ is
controlled by the finite index subgroup $\Uni(H_M)$. For these
groups, $\Uni(H_M)$ is $k$--algebraic if and only if $H_M$ is defined
over $k$.

In the real setting, these groups are of the form $\Ort(B_\iny)$,
where $B_\iny$ is a symmetric, positive definite bilinear form and the
form $B_\iny$ will be defined over a subfield $k \su
\R$. In the complex setting, these groups are of the form
$\Uni(H_\iny)$, where $H_\iny$ is a Hermitian form of signature
$(n-1,0)$ and $H_\iny$ will be defined over a subfield $k \su \C$. In
the quaternionic setting, these groups are of the form $\Uni(H_\iny)$,
where $H_\iny$ is a Hermitian form of signature $(n-1,0)$ and $H_\iny$
will be defined over a subalgebra $\Cal{A} \su \B{H}$. Our only
interest is when $k$ is a number field in the first two settings or
$\Cal{A}$ is a quaternion algebra defined over a number field in the
last setting. 

\subsection{Arithmetically defined subgroups of unitary affine groups}   

For an AB-group $\Ga$ modelled on $\Fr{N}_{\ell n-1}(X)$, we saw
in \S 2.6 that $\Ga$ can be conjugated into a subgroup of a unitary
affine group $U_M(n-1;X)$ for any $M(X)$. If this unitary affine group
is $k$--algebraic and $\Ga$ is contained in the $k$--points, we say
that $\Ga$ is \emph{$k$--defined}. When $\Ga$ is commensurable with
the $\Cal{O}_k$--points ($\Cal{O}_k$ is either the ring of integers of
$k$ or a maximal order in the quaternion algebra), we say that $\Ga$
is a \emph{$k$--arithmetic subgroup}. Note that if $\Ga$ is
$k$--defined, then by conjugating by a Heisenberg dilation, we can
arrange for $\Ga$ to be commensurable with a subgroup of the
$\Cal{O}_k$--points of the unitary affine group. 

\subsection{The quaternionic setting}

In the quaternionic setting, we can realize $U_M(n-1;\B{H})$ as
$\wh{k}$--algebraic subgroup of $\GL(m;\R)$, where $\wh{k}$ is the
field for which the quaternion algebra $\Cal{A}$ is defined.
For a maximal order $\Cal{O}$ in $\Cal{A}$ (see
\cite{MaclachlanReid03}), if $\Ga$ has a finite index in the
$\Cal{O}$--points of some unitary affine group $U_M(n-1;\B{H})$, when
we realize $U_M(n-1;\B{H})$ as a $\wh{k}$--algebraic group, $\Ga$
will have a finite index subgroup in the $\Cal{O}_{\wh{k}}$--points of
this group. 

In our notation, we will refer to $U_M(n-1;\B{H})$ as being
$\Cal{A}$--defined, subgroups $\Ga$ which are commensurable with
$\Uni(n-1;\Cal{O})$ for some maximal order $\Cal{O}$ as being
$\Cal{A}$--arithmetic, and homomorphisms as being $\Cal{A}$--defined.
Since when we realize $\Uni(n-1;\B{H})$ as a $\wh{k}$--algebraic group,
these definitions correspond to the standard algebraic definitions
(over the field $\wh{k}$), this is only a slight abuse of notation.


\subsection{$k$--monomorphisms of unitary affine groups into the  
isometry group} 

In this subsection, we characterize when a unitary affine group admits
a $k$--algebraic structure via embeddings into the isometry group of
$X$--hyperbolic space. 

Let $U_M(n-1;X)$ be a $k$--algebraic unitary
affine group. Then $H_{M(X)}$, the associated Hermitian form for
$M(X)$, is defined over $k$. Set $H = H_{M(X)} \op D_2$, with $H$
defined on $X^{n-1} \op X^2$ and $(X^2,D_2)$ is a $k$--defined
$X$--hyperbolic plane. Finally, let $V_\iny$ denote the
$H$--orthogonal complement in $X^{n+1}$ of a pair of $X$--linearly
independent, $k$--defined, $H$--null vectors $v$ and $v_0$ in $(X^2,D_2)$.

Let $\psi\co (X^{n-1},H_{M(X)}) \lra (V_\iny, H_{|V_\iny})$ be any isometric
isomorphism defined over $k$. This $X$--linear 
map induces a $k$--isomorphism
\[ \rho\co U_M(n-1;X) \lra MN, \] 
where $N$ and $M$ are factors in the Iwasawa decomposition induced on
$\Stab(v)$ with respect to the above pair of $H$--null vectors. 
Since both vectors are $k$--defined, it follows that $MN$ is $k$--algebraic.

As a result of this discussion, we have the following proposition.

\begin{prop}\label{P:311}
$U_M(n-1;X)$ is a $k$--algebraic group 
if and only if there exists a Hermitian form $H$ of signature $(n,1)$
defined over $k$ and a $k$--isomorphism
$\rho\co U_M(n-1;X) \lra MN <
\Isom(\Hy_X^n)$ where $\Hy_X^n$ is modelled on $H$.
\end{prop}

\subsection{A necessary and sufficient condition for arithmeticity}

In this subsection, we classify cusp cross-sections of arithmetic
hyperbolic lattices. In the previous subsection, we related
the algebraic structure of abstractly defined unitary affine groups
via embeddings into the isometry group of $X$--hyperbolic space. In
this subsection, we do the same for AB-groups.

We start with the following proposition which determines when an
AB-group is $k$--defined. 

\begin{prop}\label{P:321}
$\Ga$ is a $k$--defined AB-group modelled on $\Fr{N}_{\ell n-1}(X)$ if
and only if there exists a $k$--defined Hermitian form $H$ modelling
$X$--hyperbolic $n$--space, a subgroup $\La< \Uni(H;k)$ commensurable
with $\Uni(H;\Cal{O}_k)$, and an injection $\rho\co \Ga \lra \Stab(v)
\cap \La$ for some $k$--defined light-like vector $v$.
\end{prop} 

\begin{proof}
For the direct implication, assume $\Ga<U_M(n-1;X)$, 
for a $k$--defined unitary affine group. 
Let $\rho\co U_M(n-1;X) \lra MN < \Uni(H)$
be a\break $k$--isomorphism given by \refP{311}. This gives us a
$k$--monomorphism\break $\rho\co U_M(n-1;X) \lra \Uni(H)$
of $k$--algebraic groups. By \refL{323}, there exists $\La<
\Uni(H;k)$, commensurable with $\Uni(H;\Cal{O}_k)$ such that $\rho(\Ga) < \La$,
as asserted.

For the reverse implication, we assume the existence of $H$, $\La$,
$\rho$, and $v$. Note that for the Fitting subgroup $L$ of $\Ga$, $\rho(L)<N$,
for some nilpotent factor of an Iwasawa decomposition. Since $L$ is
Zariski dense in $N$ and consists of $k$--points, $N$ is a $k$--algebraic
subgroup. Since $\rho(\Ga)$ is virtually contained in $N$, $\rho(\Ga)<MN$,
for the compact factor $M$ of an Iwasawa decomposition $MAN$
of $\Stab(v)$. Since the group $M$ can be selected to be $k$--algebraic,
we have $\rho(\Ga)<MN$, where $MN$ is a $k$--algebraic unitary affine group,
as desired.
\end{proof}

For an AB-group $\Ga$ modelled on $\Fr{N}_{\ell n-1}(X)$, we say
that $\Ga$ is \emph{arithmetically admissible}  if there exists an 
arithmetic $X$--hyperbolic $n$--lattice $\La$ such that $\Ga$ is 
isomorphic to $\tri_v(\La)$. Altogether we have the following theorem which classifies the
arithmetically admissible AB-groups (part (a) is proved in \cite{LongReid02}). 

\begin{thm}[Cusp classification theorem]\label{T:324}
Let $\Ga$ be an AB-group modelled on $\Fr{N}_{\ell n-1}(X)$.
\begin{itemize}
\item[\rm(a)]
For $X=\R$, $\Ga$ is arithmetically admissible if and only if $\Ga$
is a $\Q$--arithmetic subgroup in $\R^{n-1} \rtimes \Ort(B_\iny)$, where 
$B_\iny$ is a $\Q$--defined, positive definite, symmetric bilinear form 
on $\R^{n-1}$. 
\item[\rm(b)]
For $X=\C$, $\Ga$ is arithmetically admissible if and only if $\Ga$
is a $k$--arithmetic subgroup in a unitary affine group for some
imaginary\break quadratic number field $k$.
\item[\rm(c)]
For $X=\B{H}$, $\Ga$ is arithmetically admissible if and only if
$\Ga$ is a $\Cal{A}$--arithmetic subgroup in a unitary affine group,
for some quaternion algebra $\Cal{A}$ with Hilbert symbol
$\pr{\frac{-a,-b}{\Q}}$ for $a,b \in \N$.   
\end{itemize}
\end{thm}

\begin{proof}
The direct implication is immediate in all three case. For the
converse, assume that $\Ga$ is a $k$--arithmetic subgroup in a 
unitary affine group, where $k$ is as above. By \refP{321}, there
exists a $k$--defined Hermitian form $H$ modelling 
$X$--hyperbolic $n$--space, a subgroup $\La< \Uni(H;k)$ commensurable
with $\Uni(H;\Cal{O}_k)$, and an injection $\rho\co \Ga \lra \Stab(v)
\cap \La$ for some $k$--defined light-like vector $v$.
$\rho(\Ga)$ must be a finite index subgroup of $\tri_v(\La)$ and by
\refT{Arithmetic}, $\La$ is an arithmetic subgroup. In this
injection we cannot ensure that $\rho(\Ga) = \tri_v(\La)$. As
$\La$ is an arithmetic subgroup in the $k$--algebraic group $\Uni(H)$, by
\refT{Sep}, we can find a finite index subgroup $\Pi < \La$ such that 
$\rho(\Ga) = \tri_v(\Pi)$. Specifically, select a complete set of
coset representatives for $\tri_v(\La)/\rho(\Ga)$, say
$\al_1,\dots,\al_r$. By \refT{Sep}, there exists a finite index
subgroup $\Pi$ of $\La$ such that $\rho(\Ga) < \Pi$ and for each
$j=1,\dots,r$, $\al_j \notin \Pi$. It then follows that $\tri_v(\Pi) =
\rho(\Ga)$, as desired.
\end{proof}

\begin{rem}
Using the Bieberbach theorems, we can easily see from (a) that every
Bieberbach group is arithmetically admissible. Altogether, this yields
a slightly simpler proof of the main result in \cite{LongReid02}.
\end{rem}

For an AB-group $\Ga$ modelled on $\Fr{N}_{2n-1}$, we say that the
holonomy group $\te$ of $\Ga$ is \emph{complex} if $\te \su
\Uni(H_{M(X)})<M(X)$ for the holonomy representation. Otherwise, we
say that $\te$ is \emph{anticomplex}. We have the following
alternative characterization based on the structure of the holonomy 
representation.  

\begin{cor}[Holonomy theorem]\label{C:Holonomy}
Let $\Ga$ an AB-group modelled on $\Fr{N}_{2n-1}$ with complex
holonomy. Then $\Ga$ is arithmetically admissible if and only if the
holonomy representation $\vp$ is conjugate to a representation
\[ \rho\co \te \lra \GL\pr{n-1;k} \]
for some imaginary quadratic number field.
\end{cor}

\begin{proof}
If $\Ga$ is arithmetically admissible, then from \refT{324}, there exists
a $k$--defined unitary affine group $U_M(n-1;k)$ such that $\Ga$ is
conjugate into $U_M(n-1;k)$ and commensurable with $U_M(n-1;\Cal{O}_k)$,
for some imaginary quadratic number field $k$. This yields an injective
homomorphism
\[ \rho\co \te \lra M(k). \]
Since $\te$ complex, $\rho(\te) < \Uni(H_{M(X)};k)$, which is a subgroup of
$\GL(n-1;k)$, as desired.

For the converse, assume that the holonomy representation of $\te$ maps
into $\GL(n-1;k)$, for some imaginary quadratic number field $k$.
By taking the $\te$--average of any $k$--defined Hermitian form, we see
that this representation is contained in a $k$--defined unitary group
$\Uni(H_{M(X)};k)$. Using this representation and a presentation for
$\Ga$, we get a system of linear homogenous equations with coefficients in $k$. Since
$\rho$ is conjugate to the holonomy representation, by the generalized
Bieberbach theorems, this system has a solution which yields a faithful
representation into $\Fr{N}_{2n-1}(k) \rtimes \Uni(H_{M(X)};k)$. By conjugating
by a Heisenberg dilation to ensure that the Fitting subgroup consists
of $k$--integral entries, we see that $\Ga$ is $k$--arithmetic. Therefore,
by \refT{324}, $\Ga$ is arithmetically admissible.
\end{proof}


\section{Families of examples}

In this section, we give examples which show that the
characterization of arithmetic admissibility is nontrivial. These
examples constitute a proof of \refT{Ob}.


\subsection{Prime order holonomy}

Let $\Ga$ be an AB-group modelled on $\Fr{N}_{2n-1}$ with
cyclic order $p$ holonomy, $C_p$, where $p$ is an 
odd prime. Note that this holonomy is necessarily complex
and so acts trivially on the center of the Fitting subgroup.
By taking the quotient of $\Ga$ by its center, we get an
$(2n-2)$--dimensional Bieberbach group with $C_p$--holonomy.
By the Bieberbach theorems, there exists a faithful
representation of $C_p$ into $\GL(2n-2;\Z)$. This can occur
only when $p-1\leq 2n-2$. That such AB-groups exist in dimension
$2n-1$ can be shown by explicit construction.

The following proposition shows that there are infinitely many
AB-groups modelled on $\Fr{N}_{2n-1}$ for infinitely many $n$
which are not arithmetic admissible. 

\begin{prop}\label{P:411}
Let $\Ga_p$ denote an AB-groups modelled on
$\Fr{N}_{2n-1}$ with holonomy $C_p$ and $2(n - 1) = p-1$.
If $\Ga_p$ is arithmetically admissible, then $p \equiv
3 \mod 4$.
\end{prop}

\begin{proof}
If $\Ga_p$ is arithmetically admissible, by \refC{Holonomy}, there
exists a faithful representation ($k$ is an imaginary quadratic number
field) 
\[ \rho\co C_p \lra \GL\pr{\frac{p-1}{2};k}. \] 
Let $k_\rho$ denote the field generated by the traces of $\rho(\xi)$
for a generator\break $\xi \in C_p$ and note that $k_\rho \su
k$. The representation $\rho$ is conjugate to one which decomposes
into a direct sum of characters $\chi_j\co C_p \lra \C^\times$ (see
\cite{CurtisReiner88} or \cite{Serre77}). Each of these characters
$\chi_j$ is of the form $\chi_j(\xi) = \zeta_p^{n_j}$. Therefore 
\[ \tr(\rho(\xi)) = \sum_{j=1}^{\frac{p-1}{2}} \zeta_p^{n_j}. \]
Since $\rho$ is faithful, for some $j$, $n_j \ne 0 \mod p$. By
considering the cyclotomic polynomial $\Phi_p(x)$, we deduce that
$\tr(\rho(\xi)) \notin \Q$ and so $k_\rho$ is a nontrivial extension
of $\Q$. On the other hand, from the decomposition above,
$k_\rho \su \Q(\zeta_p)$. Since $[k:\Q]=2$, it must be that
$k=k_\rho$. Hence $\Q(\zeta_p)$ contains an imaginary quadratic
extension of $\Q$. By quadratic reciprocity, this can happen if and
only if $p \equiv 3 \mod 4$.
\end{proof}

It is worth noting for $p \equiv 1 \mod 4$, \refC{Holonomy} can be used
to show that such AB-groups are arithmetically admissible and the
field $k$ is the unique imaginary quadratic number field in
$\Q(\zeta_p)$. Moreover, the holonomy generator acts by the matrix
$\res_{\Q(\zeta_p)/k}(\zeta_p)$, where $\res$ denotes the restriction
of scalar operation.  

\begin{rem}
Note when $p >5$, we get an obstruction without appealing to
\refT{324}. In this case, we have an injection 
\[ \rho\co C_p \lra \Uni(H;k) < \GL\pr{\frac{p+1}{2};k}. \]
So long as $(p+1)/2 < p-1$, this implies $p \equiv 3 \mod 4$.
\end{rem}

In the quaternionic setting, take $\Ga$ modelled on
$\Fr{N}_{7}(\B{H})$ with $C_5$--holonomy where the action
of $C_5$ on the center of $\Fr{N}_7(\B{H})$ is trivial. If $\Ga$ is
arithmetically admissible, we obtain an injection of $C_5$ into
$\Cal{A}$, a definite quaternion algebra over $\Q$. However,
this is impossible since such $\Cal{A}$ do not contain elements of
order five (see \cite{Vigneras80}). Thus $\Ga$ cannot be
arithmetically admissible. This yields: 

\begin{prop}
There exists infinitely many nonarithmetically admissible AB-groups
modelled on $\Fr{N}_7(\B{H})$.
\end{prop}

In the next subsection, by taking products of this example with
quaternionic nil-tori, we can construct infinitely many examples
in $\Fr{N}_{4n+3}(\B{H})$ for $4n+3>7$.  


\subsection{Central products}

In the complex setting, we get arithmetically inadmissible examples
in every dimension greater than two by constructing new AB-groups from a pair of
AB-groups. The details are as follows.

Let $\Ga_1$ and $\Ga_2$ denote a pair of AB-groups modelled on $\Fr{N}_{2n_1-1}$
and $\Fr{N}_{2n_2-1}$, respectively. We define $\Ga = \Ga_1 \times_{c} \Ga_2$ to 
be the group $\Ga_1 \times \Ga_2 / N_{\Ga_1\times \Ga_2}(\innp{c_1c_2^{-1}})$,
where $c_j$ is the generator of the center of the Fitting subgroup
$L_j$ of $\Ga_j$ and $N_{\Ga_1\times \Ga_2}(c_1c_2^{-1})$ denotes the
normal closure  of $\innp{c_1c_2^{-1}}$ in $\Ga_1\times\Ga_2$. We call
this the \emph{central product} of $\Ga_1$ and $\Ga_2$.

\begin{lemma}
If $\te_1$ and $\te_2$, the respective holonomy groups, are either both complex
or both anticomplex, then $\Ga_1\times_c \Ga_2$ is an AB-group modelled on 
$\Fr{N}_{2(n_1+n_2)-1}$.
\end{lemma}

\begin{proof}
As we only make use of the case when the holonomy groups are complex,
we leave the anticomplex case for the reader. To begin, we have
natural inclusions of $\Fr{N}_{2n_1-1}$ and $\Fr{N}_{2n_2-1}$ into
$\Fr{N}_{2(n_1+n_2)-1}$ induced by
\[ \rho_j\co \C^{n_j} \lra \C^{n_1} \op \C^{n_2}. \]
In fact, this yields inclusions of $\C^{n_j} \rtimes \GL(n_j;\C)$ into
$\C^{n_1+n_2} \rtimes \GL(n_1+n_2;\C)$. As a result, we have an
injective homomorphism of $\Ga_1 \times_c \Ga_2$ onto an AB-group in
$\Aff(\Fr{N}_{2(n_1+n_2)-1})$. By carefully selecting the maps
$\rho_j$, we can ensure that the induced maps agrees on the center of
the Fitting subgroups of $\Ga_1$ and $\Ga_2$. With this selection, the
induced map becomes an isomorphism of $\Ga_1 \times_c \Ga_2$ with an
AB-group in $\Fr{N}_{2(n_1+n_2)-1}$, as desired.
\end{proof}

\begin{rem}
In the complex holonomy case, the group $\innp{c_1c_2^{-1}}$ is
normal.
\end{rem} 

Using central products we can construct many arithmetically
inadmissible AB-groups. We summarize this in the following theorem
which proves \refT{Ob} (a).

\begin{thm}[Central product theorem]\label{T:ConGlue}
Let $\Ga_1$ and $\Ga_2$ be AB-groups modelled on $\Fr{N}_{2n_1-1}$ and
$\Fr{N}_{2n_2-1}$, defined over $k_1$ and $k_2$, respectively.  Assume
$\Ga_1$ and $\Ga_2$ are both complex or anticomplex.
\begin{itemize}
\item[\rm(a)]
$\Ga_1 \times_c \Ga_2$ is $k_1k_2$--defined.
\item[\rm(b)]
$\Ga_1 \times_c \Ga_2$ is $k$--arithmetic if and only if both $\Ga_1$
and $\Ga_2$ are $k$--arithmetic.
\item[\rm(c)]
There exist arithmetically inadmissible AB-groups modelled on
$\Fr{N}_{2n-1}$ for all $n \geq 3$.
\end{itemize}  
\end{thm}

\begin{proof}
For (a), let $M_1(X)$ and $M_2(X)$ denote the maximal compact groups\break
defined over $k_1$ and $k_2$ and for which $\Ga_j$ injects into the
$k_j$--points of\break $\Fr{N}_{2n_j-1}(X) \rtimes M_j(X)$. Associated to
each of these maximal compact 
groups is a $k_j$--defined Hermitian form $H_j$. Let $H = H_1 \op H_2$
and $M(X)$ be the associated maximal compact subgroup. $M(X)$ is
$k_1k_2$--define and we have an injection (into the
$k_1k_2$--points) $\rho\co \Ga_1\times_c \Ga_2 \lra
\Fr{N}_{2(n_1+n_2)-1}(X) \rtimes M(X)$. Thus $\Ga_1\times_c\Ga_2$ is
$k_1k_2$--defined.  

For (b), if $\Ga=\Ga_1 \times_c \Ga_2$ is $k$--arithmetic, then $\Ga$
is isomorphic to a maximal peripheral subgroup of an arithmetic
lattice $\La$ in $\Uni(H)$, where $H$ is a signature $(n_1+n_2,1)$
Hermitian form defined over an imaginary quadratic number field
$k$. It must be that $\Ga_j$ injects into a subgroup $\La_j$ of
$\La$ which is maximal with respect to stabilizing a $k$--defined
complex subspace $\C^{n_j,1}$. As $\La_j$ is an arithmetic lattice in
a smaller isometry group (whose model form is the restriction of $H$ to
the complex subspace $\C^{n_j,1}$), this shows that
$\Ga_j$ is $k$--arithmetic for $j=1,2$. The reverse implication
follows immediately from (a).

For (c), by \refP{411}, there exists an arithmetically inadmissible
AB-group modelled on $\Fr{N}_5$. To obtain examples in higher
dimensions, we take central products of this example with other
AB-groups and apply (b).
\end{proof}

The quaternionic setting can be handled similarly. We construct
examples in every dimension by taking central products with the
inadmissible example $\Ga$ in $\Fr{N}_7(\B{H})$ given above.

\begin{cor}\label{C:64}
There exist arithmetically inadmissible AB-groups modelled on
$\Fr{N}_{4n-1}(\B{H})$ for all $n \geq 2$.
\end{cor}

\begin{rem}
\refT{Ob} (b) shows that a quaternionic version of \refCo{11} (a) is
false. Namely, there exist almost flat manifolds modelled on the
quaternionic Heisenberg group which cannot arise as the a cusp
cross-section of a 1-cusped quaternionic hyperbolic manifold. The
author has been informed by Walter Neumann and Alan Reid that nil
3--manifold groups exist which cannot arise as a cusp cross-section of
a 1--cusped complex hyperbolic 2--manifold. In fact, Yoshinobu
Kamishima \cite{Kamishima04} proves that no nil 3--manifold group with
nontrivial holonomy can arise as a cusp cross-section of a 1--cusped
complex hyperbolic 2--manifold.
\end{rem}


\section{Analysis in low complex dimensions}

In this section, we work out the specific details of \refT{324} for
closed nil 3--manifold groups. In particular, we prove
\refC{Nil3}. Before we undertake this task, it is worth noting that
the validity of \refC{Nil3} follows almost at once from
\refC{Holonomy}. The only ingredient needed is the classification of
nil 3--manifold groups given below. However, in this section we
explicitly construction the desired $k$--integral representations for
such groups to better illustrate the proof of \refT{324}.


\subsection{Some preliminary material}

In this subsection, we give an alternative (more traditional) description
of nil 3--manifolds and list the classification of closed nil 3--manifold
groups. See \cite{Scott83} for a reference on this material. 

The 3--dimensional Heisenberg group also has a definition as
\[ \Fr{N}_3 = \set{\begin{pmatrix} 1 & x & t \\ 0 & 1 & y \\ 0 & 0 & 1
\end{pmatrix} ~:~ x,y,t \in \R}. \]
We will identify $S^1$ with the rotations in the $xy$--plane. An
\emph{orientable nil 3--manifold} is a manifold of the form
$\Fr{N}_3/\Ga$, where $\Ga$ is a discrete subgroup of $\Fr{N}_3
\rtimes S^1$ which acts freely. As we
will have need for this in the sequel, we must also consider an
orientation reversing involution given by
\[ \wt{\iota}\begin{pmatrix} 1 & 0 & 0 \\ 0 & 1 & y \\ 0 & 0 & 1
\end{pmatrix} = \begin{pmatrix} 1 & 0 & 0 \\ 0 & 1 & -y \\ 0 & 0 & 1
\end{pmatrix}, \quad \wt{\iota} \begin{pmatrix} 1 & x & 0 \\ 0 & 1 & 0
\\ 0 & 0 & 1 \end{pmatrix} = \begin{pmatrix} 1 & x & 0 \\ 0 & 1 & 0 \\
0 & 0 & 1 \end{pmatrix}. \]
As automorphisms of lattices of $\Fr{N}_3$ uniquely determine
automorphisms of $\Fr{N}_3$ by Mal'cev rigidity, this
determines a continuous isomorphism of $\Fr{N}_3$.  A \emph{nil
3--manifold} is a manifold of the form $\Fr{N}_3/\Ga$, where $\Ga$ is a
discrete subgroup of $\innp{\Fr{N}_3\rtimes S^1,\wt{\iota}}$ which acts
freely.

\begin{rem}
This terminology follows \cite{Scott83} and
\cite{Thurston97}. However, in our terminology which follows
\cite{Dekimpe96}, such manifolds should be called infranil
3--manifolds. We adhere to the first in this section and hope no
confusion arises from this. 
\end{rem}    

In \S 2, we gave a general definition of the
Heisenberg group. That these two definitions coincide follows from the
fact that both groups are connected, simply connected 2--step nilpotent
Lie groups of dimension three and there is a unique such Lie group (up to
Lie isomorphism) that has these properties. More generally, any simply
connected, connected, 2--step nilpotent group with 1--dimensional center
is uniquely determined (see \cite{McReynolds04C}) and is Lie isomorphic to the
the Heisenberg group which we defined in \S 2.

In this form, we identify $S^1$ with $\Uni(1)$, where $\Uni(1)$ acts
(as above) by 
\[ U(z,t) = (Uz,t). \]
We also identify $\iota$ and $\wt{\iota}$, where $\iota$ is the isometry
induced by complex conjugation. When the nil 3--manifold is
closed, $\Fr{N}_3 \cap \Ga$ must be finite index in $\Ga$ and
$\Fr{N}_3/(\Fr{N}_3\cap\Ga)$ must be compact.

The following is a complete list of closed nil 3--manifold groups (see
\cite[p.\ 159--166]{Dekimpe96}).

\begin{smallgroup}
\begin{itemize}
\item[(1)]
\[ \langle a,b,c~:~[b,a] = c^k,~[c,a] = [c,b] = 1\rangle. \]
with $k \in \N$.
\item[(2)]
\begin{align*}
\langle a,b,c,\al~:~&[b,a]= c^k,~[c,a] = [c,b] = [\al,c] = 1,~\al a =
a^{-1} \al \\ & ~\al b = b^{-1}\al,~\al^2 = c\rangle,
\end{align*}
with $k \in 2\N$.
\item[(3)]
\begin{align*}
\langle a,b,c,\al~:~&[b,a] = c^{2k},~
[c,a] = [c,b] = [a,\al] =  1, ~\al c = c^{-1}\al,\\
& ~\al b = b^{-1}\al c^{-k},~\al^2 = a\rangle
\end{align*}
with $k \in \N$.
\item[(4)]
\begin{align*}
 \langle a,b,c,\al,\be~:~& [b,a]=c^{2k},~
[c,a]=[c,b]=[c,\al]=[a,\be]=1,\\ &~\be c= c^{-1}\be,~ \al a = a^{-1}
\al c^k, \al b = b^{-1}\al c^{-k},\\ &~ \al^2 = c,~ \be^2 = a, ~ \be b
= b^{-1}\be c^{-k},~ \al\be = a^{-1}b^{-1}\be \al c^{-k-1}\rangle, 
\end{align*} 
with $k \in 2\N$.
\item[(5)]
\begin{align*}
\langle a,b,c,\al~:~&[b,a] = c^k,~[c,a]  = [c,b] = [c,\al] = 1,~ \al a  =
b\al,\\ &~\al b  = a^{-1}\al,~\al^4  = c^p\rangle,
\end{align*}
with $k \in 2\N$ and $p=1$ or $k\in 4\N$ and $p=3$.
\item[(6)]
\begin{align*}
\langle a,b,c,\al~:~&[b,a] = c^k,~
[c,a] = [c,b] = [c,\al] = 1,~
\al a = b\al c^{k_1},\\
&~\al b = a^{-1}b^{-1}\al,~
\al^3 = c^{k_2}\rangle
\end{align*}
with $k>0$ and
\[ k \equiv 0 \mod 3, \quad k_1=0,~k_2 = 1, \text{ or } k \equiv 0 \mod
3, \quad k_1=0,~k_2 = 2, \]  
or
\[ k \equiv 1,2 \mod 3, \quad k_1=1,~k_2=1. \]
\item[(7)]
\begin{align*}
\langle a,b,c,\al~:~&[b,a] = c^k,~ [c,a] = [c,b] = [c,\al] = 1,~\al a =
ab\al,\\
&~\al b = a^{-1}\al,~\al^6 = c^{k_1}\rangle,
\end{align*}
with $k>0$ and
\[ k \equiv 0 \mod 6, \quad k_1 = 1,
\text{ or }
k \equiv 4 \mod 6, \quad k_1 = 1, \]
or
\[ k \equiv 0 \mod 6, \quad k_1 = 5,
\text{ or } 
k \equiv 2 \mod 6, \quad k_1 = 5. \]
\end{itemize}
\end{smallgroup}


\subsection{Nil 3--manifolds}

Let $M$ be a closed nil 3--manifold with $\pi_1(M)=\Ga$. By
\refT{GenBieberbach}, in order to prove \refC{Nil3}, it suffices to
show that $\Ga \cong \tri(\La)$, where $\La$ is an arithmetic
lattice in $\Isom(\Hy_\C^2)$. In fact, by \refT{Cusp} it suffices
to construct an injective homomorphism $\vp\co \Ga \lra \tri(\La)$.                      

Let $\Fr{N}(3) = \Fr{N}_3\rtimes \Uni(1)$ and $\iota$ be the
isometry of $\Hy_\C^2$ induced by conjugation. For a subring $R \su
\C$, we define $\Fr{N}_3(R) = R \times \Ima R$ with the induced group
operation and set $\Fr{N}(3,R) = \Fr{N}_3(R) \rtimes \Uni(1;R)$. 

For the statement, let $\zeta_3$ be a primitive third root of unity,
say$$\zeta_3=-1/2+\sqrt{-3}/2.$$  

\begin{thm}\label{T:22}
Let $M$ be a closed nil 3--manifold $M$ and $\Ga =\pi_1(M)$. Then there
exists a faithful representation $\vp\co \Ga \lra
\innp{\Fr{N}(3,\Cal{O}_k),\iota}$ with $k=\Q(i)$ or $\Q(\zeta_3)$.
\end{thm}

\begin{proof}
We begin by summarizing the strategy of the proof, which depends
heavily on the list of presentations for the
fundamental group of a closed nil 3--manifold in previous subsection.  
The idea is to show that an injective homomorphism on the Fitting
subgroup $\innp{a,b,c}$ can be promoted to the full 3--manifold group 
(see \cite[Thm 3.1.3]{Dekimpe96}). To get a representation with the coefficients 
in $\Z[i]$ or $\Z[\zeta_3]$, we are reduced to solving some simple equations. 
The details are as follows. 

In the lemma below, let $p_1\co \Fr{N}_3 \lra \C$ be projection onto
the first factor. 

\begin{lemma}\label{L:311}$\phantom{9}$
\begin{itemize}
\item[\rm(a)]
Let $a$ and $b$ be as above and $\rho\co \innp{a,b,c} \lra \Fr{N}_3$ 
be a homomorphism. If $p_1(\rho(a))$ and $p_1(\rho(b))$ are
$\Z$--linearly independent and $c\notin \ker \rho$, then $\rho$ is
injective. 
\item[\rm(b)]
Let $\rho\co \Ga \lra \innp{\Fr{N}(3),\iota}$
be a homomorphism. If $\rho^{-1}(\rho(\innp{a,b,c})) = \innp{a,b,c}$,
and $\rho_{|\innp{a,b,c}}$ is an injective homomorphism, then $\rho$
is an injective homomorphism. 
\end{itemize}
\end{lemma}

\begin{proof}
For (a), let $w\in \ker \rho$ and write $w = a^{n_1}b^{n_2}c^{n_3}$.
Also, let\break $\rho(a) = (v_1,t_a)$, $\rho(b)=(v_2,t_b)$. Since
$[a,b]=c^k$, it must be that $\rho(c)=(0,s)$ as
$[\Fr{N}_3,\Fr{N}_3]=\set{(0,t)~:~t \in \R}$. With this said, we
see that 
\[ \rho(w) =
(n_1v_1+n_2v_2,n_1t_a+n_2t_b+2\Ima\innp{n_1v_1,n_2v_2}+n_3s). \] 
Since $w \in\ker \rho$, $n_1v_1 + n_2v_2 = 0$. The
$\Z$--linear independence of $v_1$ and $v_2$ implies
$n_1=n_2=0$. Therefore $n_3s=0$, and so $n_3=0$, as $s \ne 0$.

For (b), let $w \in \ker \rho$ and write
\[ w = a^{n_1}b^{n_2}c^{n_3}\al^{s_1}\be^{s_2}, \quad
n_1,n_2,n_3,s_1,s_2\in \Z. \]
Since $\rho(w)=1$, using the above form for $w$, we have
\[ \rho(a^{n_1}b^{n_2}c^{n_3}) = \rho(\al^{-s_1}\be^{-s_2}). \]
Since $\rho^{-1}(\rho(\innp{a,b,c}))=\innp{a,b,c}$, it must be that
$\al^{s_1}\be^{s_2} \in \innp{a,b,c}$. Thus\break $w \in
\innp{a,b,c}$. Since $\rho_{|\innp{a,b,c}}$ is one-to-one, $w=1$, as 
desired.
\end{proof}

\begin{rem}
Notice that for (b), we only require that if $\rho(\al^{s_1}\be^{s_2})
\in \rho(\innp{a,b,c})$ then $\al^{s_1}\be^{s_2} \in
\innp{a,b,c}$. Moreover, we only need to check for $s_1 \in
\set{1,\dots,k_\al}$ and $s_2 \in \set{1,\dots,k_\be}$ where $k_\eps$
is the first integer such that $\eps^{k_\eps}\in\innp{a,b,c}$, $\eps =
\al$ or $\be$.
\end{rem}

Let $L=\innp{a,b,c}$ and define two homomorphisms
$\vp_3,\vp_4\co L \lra \Fr{N}_3$
by\break $\vp_j(a) = (1,0)$ and $\vp_j(b) = (\zeta_j,0)$.
This determines $c$, since some power of $c$ is a commutator of $a$ and
$b$. By \refL{311} (a), both maps are injective homomorphisms.

We extend this to $\Ga$ by declaring
\begin{equation}\label{E:311}
\vp_j(\al) = (z_1,t_1,\eta_1), \quad \vp_j(\be) = (z_2,t_2,\eta_2),
\end{equation}
where $z_1,z_2\in\C$, $t_1,t_2\in\R$, and
$\eta_1,\eta_2\in\innp{\Uni(1),\iota}$. 

To solve \refE{311}, we simply use the presentations above to ensure
that this yields a homomorphism. By applying \refL{311} (b),
one can see that these solutions yield injective homomorphisms. For
clarity, we solve the equations for the second family (2) and give a
list of the equations and solutions for the seventh family (7). We
relegate the rest of the solutions to the appendix.

The second family has presentation
\begin{smallgroup}
\[ \innp{a,b,c,\al~:~ [b,a]=c^k,~[c,a]=[c,b]=[c,\al]=1,~\al a = a^{-1}
\al, \al b = b^{-1} \al, \al^2 = c} \]
\end{smallgroup} 
with $k \in 2\N$. For this family we take the map $\vp_4$. First,
consider the relation $[b,a]=c^k$. Then
\begin{align*}
[\vp_4(b),\vp_4(a)] &= [(i,0,1),(1,0,1)] \\
&= (i+1,2\Ima\innp{i,1},1)(-i-1,2\Ima\innp{-i,-1},1) \\
&= (0,4,1).
\end{align*}
Since $[\vp_4(b),\vp_4(b)]=\vp_4(c)^k$, it follows that
$\vp_4(c)=(0,4/k,1)$. 

Next, consider the relation $\al^2 = c$. We have
\begin{align*}
(z_1,t_1,\eta_1)(z_1,t_1,\eta_1) &= (z_1 + \eta_1z_1,2t_1 +
2\Ima\innp{z_1,\eta_1z_1},\eta_1^2) \\
&= (0,4/k,1).
\end{align*}
In particular, $\eta_1^2=1$. If $\eta_1=1$, then the above injection
would yield an isomorphism between a group in the first family with a
group in the second family. This is impossible, therefore
$\eta_1=-1$. By considering the first coordinate equation with
$\eta_1=-1$, we get no information. The second coordinate equation is
$2t_1 = 4/k$, therefore $t_1=2/k$. One can now check that $[c,\al]=1$,
regardless of $z_1$. 

Now, we take the relation $\al a = a^{-1} \al$. We have
\begin{align*}
(z_1,2/k,-1)(1,0,1) &= (z_1 - 1, 2/k - 2\Ima z_1, -1).
\end{align*}
On the other hand,
\begin{align*}
(-1,0,1)(z_1,2/k,-1) &= (-1+z_1,2/k + 2\Ima z_1,-1).
\end{align*}
The first and last coordinate equations yield no information, while
the second coordinate equation yields $4\Ima z_1=0$. Hence $\Ima
z_1=0$. 

Lastly, we have the relation $\al b = b^{-1}\al$. We have
\begin{align*}
(z_1,2/k,-1)(i,0,1) &= (z_1 - i,2/k + 2\RE z_1,-1).
\end{align*}
On the other hand,
\begin{align*}
(-i,0,1)(z_1,2/k,-1) &= (-i+z_1,2/k - 2\RE z_1,-1).
\end{align*}
As above, the first and second coordinates yields no information,
while the second coordinate implies that $\RE z_1=0$.

Hence, we deduce from the above computations, the desired homomorphism
$\rho$ is defined by
\begin{align*}
\rho(a) &= (1,0,1), & \rho(b) = (i,0,1) & \\
\rho(c) &= (0,4/k,1) & \rho(\al) = (0,2/k,-1). & 
\end{align*}

Notice that these solutions are in $\Q(i)$ and not $\Z[i]$. This is
rectified by conjugating the above representation by a 
dilation of $2k$. This dilation is linear on the first factor,
quadratic on the second factor, and trivial on the third factor. The
resulting faithful representation is
\begin{align*}
\rho(a) &= (2k,0,1), & \rho(b) = (2ki,0,1) & \\
\rho(c) &= (0,16k,1) & \rho(\al) = (0,8k,-1). & 
\end{align*}
That this is faithful follows from \refL{311} (b) once the conditions
of this lemma have been verified. The injectivity of
$\rho_{|\innp{a,b,c}}=\vp_4$ follows from\break \refL{311} (a). To check
that $\rho(\rho^{-1}(\innp{a,b,c}))=\innp{a,b,c}$, by the remark
proceeding the proof of \refL{311}, it suffices to show that $\al
\notin \rho(\rho^{-1}(\innp{a,b,c}))$. The validity of this is
obvious. 

 
For the seventh family (7), we have the presentation
\begin{smallgroup}
\[ \langle a,b,c,\al~:~[b,a]=c^k,~[c,a]=[c,b]=[c,\al]=1,~\al a =
ab\al,~\al b=a^{-1}\al,~\al^6=c^{k_1}\rangle \]
with
\begin{align*}
k \equiv 0 \mod 6,~k_1 =1, & \text{ or } k \equiv 4 \mod 6,~k_1=1,
\text{ or }\\
k \equiv 0 \mod 6,~k_1 =5, & \text{ or } k \equiv 2 \mod 6,~k_1=5.
\end{align*}
\end{smallgroup}
We take $\vp_3$ in this case. By considering all the relations, we get
\begin{align*}
(0,4/k,1) &= (0,s,1)\\
(z_1+\eta_1,t_1+2\Ima\innp{z,\eta_1},\eta_1) &= (1+\zeta_3 + z_1, t_1
+  \\
& \quad \quad \quad 2\Ima\innp{1+\zeta_3,z_1} + 2\Ima\innp{1,\zeta_3}
+ t_1,\eta_1) \\
(z_1 + \eta_1\zeta_3,t_1+2\Ima\innp{z_1,\eta_1\zeta_3},\eta_1) &=
(-1+z_1,t_1+2\Ima\innp{-1,z_1},\eta_1).
\end{align*}
We omit the last relation $\al^6 =c^{k_1}$ as it is quite long.
Note the commutator relations (aside from
$[b,a]=c^k$) are all trivially satisfied. Solving these equations and
conjugating by a dilation of $12k$ to get the coefficients in
$\Z[\zeta_3]$, we have
\begin{align*}
\vp_3(a) &= (12k,0,1), & \vp_3(b)= (12k\zeta_3,0,1), \\ \vp_3(c) &=
(0,288k\sqrt{3},1),
& \vp_3(\al)= (-6k,12k(4k_1+3k)\sqrt{3},\zeta_6). &
\end{align*}
That this is faithful follows from \refL{311} and the remark following
the proof.  
 
The other remaining families are handled similarly.
\end{proof}

The proof is completed by following the injection $\rho$ from
\refT{22} with the injection
\[ \psi\co \Fr{N}_3 \rtimes \innp{\Uni(1),\iota} \lra \Uni(2,1) \]
given by
\[ (\xi,t,\iota^\eps U) \lmto \begin{pmatrix} 1 & \xi & \xi \\ -\ol{\xi} &
1-\frac{1}{2}(\norm{\xi}^2 -it) & -\frac{1}{2}(\norm{\xi}^2-it) \\
\ol{\xi} & \frac{1}{2}(\norm{\xi}^2-it) & 1+\frac{1}{2}(\norm{\xi}^2-it)
\end{pmatrix} \begin{pmatrix} U & 0 & 0 \\ 0 & 1 & 0 \\ 0 & 0 & 1
\end{pmatrix} \iota^\eps. \]
For completeness, we give all the linear representations in the
appendix.

\begin{rem}
\refT{22} can be strengthened somewhat. Specifically, we can prove
that there is a dense set of almost flat structures on any compact nil
3--manifold which can be realized in cusp cross-sections of arithmetic
complex hyperbolic 2--orbifolds. We refer the reader for
\cite{McReynolds04C} for more on this.
\end{rem} 


\subsection{Nil 5--manifolds}

Using the list of isomorphism types of holonomy groups for nil
5--manifolds given in \cite{DekimpeEick02}, we can carry out the same analysis. 

\begin{prop}
The only complex holonomy groups which yield arithmetically
inadmissible groups are $C_5$, $C_{10}$, $C_{12}$ and $C_{24}$.
\end{prop}

To see this, by \refC{Holonomy}, it suffices to check the field of
definition for all the distinct representations of the holonomy group
in $\GL(2;\C)$. For most of the groups, 
every representation will be conjugate to one defined over
an imaginary quadratic number field. Both $C_{12}$ and $C_{24}$ can
arise via central products of nil 3--manifold groups for which
\refT{ConGlue} can be applied to show that the resulting AB-groups are
arithmetically inadmissible. For $C_5$--holonomy, we can apply \refP{411},
which in turn yields the result for $C_{10}$--holonomy, since any
arithmetic representation for an AB-group with $C_{10}$--holonomy would yield
one for an AB-group with $C_5$--holonomy.

\section{Final comment}

As discussed in the remark following \refC{64}, one cannot expect to obtain
manifolds of finite volume with just a single cusp in the above
construction. However, a natural question is whether or not one can 
guarantee that $\Hy_\C^2/\Pi$
be constructed to be a manifold. For instance, in the second family (2),
the only congruence subgroup which contains $\Ga$ is the level two
congruence subgroup. This is seen to have torsion though, as
\[ \begin{pmatrix} -1 & 0 & 0 \\ 0 & 1 & 0 \\ 0 & 0 & 1 \end{pmatrix}
\]
is an element of this subgroup. However by results of Borel \cite{Borel63},
we can guarantee that the first and third families map into torsion
free groups. We can simply conjugate by a dilation to ensure that the
representation maps into a congruence subgroup which is torsion
free. \cite{Borel63} is needed in showing that the congruence subgroup is
torsion free.


\section{Appendix: solutions to \refE{311}}

Here, we include the matrices found by solving \refE{311} in \S 7.
For all possible faithful representations, see \cite{McReynolds04C}.

\begin{smallgroup}
(1)
\begin{align*}
a &= \begin{pmatrix} 1 & 2k & 2k \\ -2k & 1- 2k^2 &
-2k^2 \\ 2k & 2k^2 & 1 + 2k^2 \end{pmatrix}, & b = \begin{pmatrix} 
1 & 2ki & 2ki \\ 2ki & 1-2k^2 & 2k^2 \\ -2ki & 2k^2 & 1+2k^2 \end{pmatrix} & \\
c &= \begin{pmatrix} 1 & 0 & 0 \\ 0 & 1 +8ki & 8ki \\ 0 & -8ki & 1-8ki
\end{pmatrix}.
\end{align*}
(2)
\begin{align*}
a &= \begin{pmatrix} 1 & 2k & 2k \\ -2k & 1- 2k^2 &
-2k^2 \\ 2k & 2k^2 & 1 + 2k^2 \end{pmatrix}, &
b = \begin{pmatrix} 1 & 2ki & 2ki \\ 2ki & 1-2k^2 &
2k^2 \\ -2ki & 2k^2 & 1+2k^2 \end{pmatrix}, & \\
c &= \begin{pmatrix} 1 & 0 & 0 \\ 0 & 1 +8ki & 8ki \\ 0 & -8ki & 1-8ki
\end{pmatrix}, &
\al = \begin{pmatrix} -1 & 0 & 0 \\ 0 & 1+ki & ki \\ 0 & -ki & 1-ki
\end{pmatrix}. &
\end{align*} 
(3)
\begin{align*}
a &= \begin{pmatrix} 1 & 4k & 4k \\ -4k & 1- 8k^2 &
-8k^2 \\ 4k & 8k^2 & 1 + 8k^2 \end{pmatrix}, &
b = \begin{pmatrix} 1 & 4ki & 4ki \\ 4ki & 1-8k^2 &
8k^2 \\ -4ki & 8k^2 & 1+8k^2 \end{pmatrix}, & \\
c &= \begin{pmatrix} 1 & 0 & 0 \\ 0 & 1 +16ki & 16ki \\ 0 & -16ki &
1-16ki \end{pmatrix}, &
\al = \begin{pmatrix} 1 & 2k & 2k \\ 2k & 1 - 2k^2 & -2k^2 \\ -2k &
2k^2 & 1 + 2k^2 \end{pmatrix}\iota. &
\end{align*}
(4)
\begin{align*}
a &= \begin{pmatrix} 1 & 4k & 4k \\ -4k & 1- 8k^2 &
-8k^2 \\ 4k & 8k^2 & 1 + 8k^2 \end{pmatrix}, &
b = \begin{pmatrix} 1 & 4ki & 4ki \\ 4ki & 1-8k^2 &
8k^2 \\ -4ki & 8k^2 & 1+8k^2 \end{pmatrix}, & \\
c &= \begin{pmatrix} 1 & 0 & 0 \\ 0 & 1 +16ki & 16ki \\ 0 & -16ki & 1-16ki
\end{pmatrix}, &
\al = \begin{pmatrix} -1 & 2k + 2ki & 2k+2ki \\ -2k+2ki & 1-4k^2 &
-4k^2 \\ 2k-2ki & 4k^2 & 1+4k^2 
\end{pmatrix}, & \\
\be &= \begin{pmatrix} 1 & 2k & 2k \\ -2k & 1-4k^2 & -4k^2 \\ 2k &
4k^2 & 1+4k^2 \end{pmatrix} \iota.
\end{align*}
(5)
\begin{align*}
a &= \begin{pmatrix} 1 & 2k & 2k \\ -2k & 1- 2k^2 &
-2k^2 \\ 2k & 2k^2 & 1 + 2k^2 \end{pmatrix}, &
b = \begin{pmatrix} 1 & 2ki & 2ki \\ 2ki & 1-2k^2 &
2k^2 \\ -2ki & 2k^2 & 1+2k^2 \end{pmatrix}, & \\
c &= \begin{pmatrix} 1 & 0 & 0 \\ 0 & 1 +8ki & 8ki \\ 0 & -8ki & 1-8ki
\end{pmatrix}, &
\al = \begin{pmatrix} i & 0 & 0 \\ 0 & 1+\frac{pk}{2}i &
\frac{pk}{2}i \\ 0 & -\frac{pk}{2}i & 1-\frac{pk}{2}i
\end{pmatrix}. &
\end{align*}
(6)
\begin{align*}
a &= \begin{pmatrix} 1 & 24k & 24k \\ -24k & 1 - 288k^2 & -288k^2
\\ 24k & 288k^2 & 1+ 288k^2 \end{pmatrix}, \\
b &= \begin{pmatrix} 1 & -12 + 12\sqrt{3}i & -12 + 12\sqrt{3}i \\
12 + 12\sqrt{3}i & 1 -288k^2 & -288k^2 \\ -12-12\sqrt{3}i & 288k^2 &
1+288k^2 \end{pmatrix},  \\
c &= \begin{pmatrix} 1 & 0 & 0 \\ 0 & 1 + 144k\sqrt{3}i &
144k\sqrt{3}i \\ 0 & -144k\sqrt{3}i & 1-144k\sqrt{3}i \end{pmatrix}, \\
\al &= \begin{pmatrix} 1 & \mu &
\mu \\ 6[k+2k_1]+6\sqrt{3}[k-2k_1]i &
1-\si & -\si \\ -6[k+2k_1]-6\sqrt{3}[k-2k_1]i & \si & 1 +\si
\end{pmatrix} \begin{pmatrix} \zeta_3 & 0 & 0 \\ 0 & 1 & 0 \\ 0 & 0 &
1 \end{pmatrix},
\end{align*}
where
\[ \si = \frac{1}{2}\brac{36(k+2k_1)^2 + 108(k-2k_1)^2 -
 192k\sqrt{3}(k\norm{z}^2 + 2k_2)i} \]
and
\[ \mu = -6[k+2k_1] + 6\sqrt{3}[k-2k_1]i. \]
(7)
\begin{align*}
a &= \begin{pmatrix} 1 & 12k & 12k \\ -12k & 1-144k^2 & -144k^2 \\ 12k
& 144k^2 & 1+144k^2 \end{pmatrix}, \\
b &= \begin{pmatrix} 1 & -6k + 6k\sqrt{3}i & -6k + 6k\sqrt{3}i \\ 6k +
6k\sqrt{3}i & 1-144k^2 & -144k^2 \\ -6k - 6k\sqrt{3}i & 144k^2 &
1+144k^2 \end{pmatrix}, \\
c &= \begin{pmatrix} 1 & 0 & 0 \\ 0 & 1 + 288k\sqrt{3}i &
288k\sqrt{3}i \\ 0 & -288k\sqrt{3}i & 1-288\sqrt{3}i \end{pmatrix}, \\
\al &= \begin{pmatrix} 1 & - 6k & -6k \\ 3k-3k\sqrt{3}i & 1 -\chi &
-\chi\\ -3k+3k\sqrt{3}i & \chi & 1+ \chi \end{pmatrix}\begin{pmatrix}
\zeta_6 & 0 & 0 \\ 0 & 1 & 0 \\ 0 & 0 & 1 \end{pmatrix}, 
\end{align*}
where
\[ \chi =  -36k^2 + 12k(4k_1+3k)\sqrt{3}i. \]
\end{smallgroup}



\Addresses\recd

\end{document}